\documentclass[a4paper,11pt]{amsart}
\usepackage[english]{babel}
\usepackage{amssymb}
\usepackage{enumitem}
\usepackage{longtable}
\usepackage[hidelinks]{hyperref}
\newtheorem{theorem}{Theorem}[section]
\newtheorem{proposition}[theorem]{Proposition}
\newtheorem{lemma}[theorem]{Lemma}
\newtheorem{corollary}[theorem]{Corollary}
\theoremstyle{definition}
\newtheorem{remark}[theorem]{Remark}
\newcommand{\Ric}{\operatorname{Ric}}

\renewcommand{\div}{\operatorname{div}}
\numberwithin{equation}{section}

\title[Symmetry breaking Poincar\'e--Einstein metrics]{Symmetry breaking Poincar\'e--Einstein metrics\\ and
non-extendable boundary Killing fields}
\author{Erwann Delay}
\address{Avignon Universit\'e, Laboratoire de Math\'ematiques d'Avignon (EA 2151), F-84916 Avignon, France}
\email{erwann.delay@univ-avignon.fr}
\date{November 23, 2026}
\begin{document}
\begin{abstract}
We construct Poincar\'e--Einstein metrics on $\mathbb R^2\times S^3\times S^2\times S^2$ which bifurcate from
anti-de Sitter black holes whose horizon is the product of a B\"ohm metric on $S^3\times S^2$ with a round $S^2$, and
which break the rotational symmetry of the $S^2$ factor while keeping their conformal infinity. Some Killing fields of
the conformal infinity therefore do not extend to the bulk, although $\pi_1(M,\partial M)=0$. The construction rests on
two results of independent interest: the Lichnerowicz eigenvalues of B\"ohm metrics on invariant
transverse traceless tensors are simple (with at most two exceptional values), the lowest one tends to $-\infty$ along B\"ohm's sequence with a gap to the
next one tending to infinity; and large anti-de Sitter black holes in dimension nine with an Einstein horizon have no static
$L^2$ perturbation of scalar type, which we prove with an explicit positive supersolution of the
Kodama--Ishibashi master equation.
\end{abstract}
\maketitle
\setcounter{tocdepth}{1}
\tableofcontents

\section{Introduction}\label{sec:intro}

Let $(M^{n+1},g)$ be a Poincar\'e--Einstein metric, i.e. a conformally compact metric with $\Ric(g)=-ng$, and let
$X_0$ be a conformal Killing field of its conformal infinity $(\partial M,[\gamma])$. When does $X_0$ extend to a Killing
field of $g$? Formally, the extension is always possible: there is a unique vector field $X$ asymptotic to $X_0$ solving
the Killing extension equation $(\nabla^*\nabla+n)X^\flat=0$, and $X$ is Killing if and only if the infinitesimal
Einstein deformation $\mathcal L_Xg$ vanishes. When $g$ is non-degenerate, i.e. $\Delta_L+2n$ has no $L^2$ kernel ($\Delta_L$ is the Lichnerowicz Laplacian, see section \ref{sec:notation}),
this always holds (this is classical, see e.g. \cite[Prop.~A.1]{ACD2}). The purpose of this paper is to construct (degenerate) Poincar\'e--Einstein metrics
for which the extension fails; the mechanism is a symmetry breaking bifurcation.

The starting point is the family of Riemannian Schwarzschild--anti-de Sitter black holes, also called (generalised)
Kottler metrics \cite{Kottler,Bir}, with an Einstein horizon $(B^N,g_B)$,
$\Ric_B=(N-1)g_B$, $n=N+1$:
\begin{equation}\label{bh}
g_{r_h}=f\,d\tau^2+\frac{dr^2}f+r^2g_B,\qquad f(r)=1+r^2-\frac{r_h^{N-1}(1+r_h^2)}{r^{N-1}},\qquad r\ge r_h ,
\end{equation}
$\tau\in\mathbb R/\tau_0\mathbb Z$, $\tau_0=4\pi/f'(r_h)$. They are Poincar\'e--Einstein metrics on
$M=\mathbb R^2\times B$, with conformal infinity $[d\tau^2\oplus g_B]$ on $\partial M=S^1_{\tau_0}\times B$, and
\begin{equation}\label{pi1}
\pi_1(M,\partial M)=0,
\end{equation}
since $\pi_1(S^1\times B)\to\pi_1(\mathbb R^2\times B)=\pi_1(B)$ is onto. The family is real-analytic in $r_h$
(after rescaling $\tau$), and every isometry of $g_B$ acts on it. As observed by Gibbons and Hartnoll \cite{GH},
a transverse traceless eigentensor of the Lichnerowicz Laplacian of $B$ with eigenvalue below
$\lambda_c=\frac{(N-1)(9-N)}4$ produces a zero mode of $\Delta_L+2n$ for some value of $r_h$, and B\"ohm metrics
\cite{Bo} provide horizons with very negative eigenvalues \cite{GHP,Ve}. If such a zero mode $k$ is charged under an
isometry of $B$, that is, not invariant under it ($\mathcal L_Xk\ne0$ for the corresponding Killing field $X$), one may
hope for a branch of Einstein metrics breaking this isometry. We take
$$
B=B_1\times S^2,\qquad g_B=g_{B_1}\oplus g_{S^2},\qquad N=7,\ n=8,
$$
where $g_{B_1}$ is a B\"ohm metric on $B_1:=S^3\times S^2$ and $g_{S^2}$ a round metric, both factors being scaled to
have Einstein constant $E:=N-1=6$.
Recall that B\"ohm \cite{Bo} constructed an infinite sequence $(g_i)_{i\ge1}$ of pairwise non-homothetic Einstein
metrics on $S^3\times S^2$, invariant under $G_1=SO(3)\times SO(3)$, which acts with cohomogeneity one and principal
orbits $S^2\times S^2$. Such a metric reads
$$
g_i=dt^2+a_i(t)^2g^{S^2}\oplus b_i(t)^2g^{S^2},\qquad t\in[0,T],
$$
where $g^{S^2}$ is the round metric of curvature $1$, $t$ is the arclength along the geodesics orthogonal to the
orbits, $T=T_i$ is the distance between the two singular orbits $\{t=0\}$ and $\{t=T\}$, the first sphere collapses
at both of them ($a_i(0)=a_i(T)=0$, $b_i>0$), and these metrics are symmetric, $a_i(T-t)=a_i(t)$,
$b_i(T-t)=b_i(t)$ \cite[\S3, \S6]{Bo}. Let $\rho$ be the
isometry of $B_1$ induced by $t\mapsto T-t$; it commutes with $G_1$. Let $\hat G_1:=G_1\times\langle\rho\rangle$,
$\Sigma=SO(3)$ acting on $S^2$, $K\subset\Sigma$ the rotations about the $z$-axis, and $\sigma\in\Sigma$ the rotation of
angle $\pi$ about the $x$-axis, which normalises $K$. Put
$$
H:=K\times\hat G_1\times O(2)_\tau ,
$$
$O(2)_\tau$ being generated by the translations and the reflection $\tau\mapsto-\tau$. All these groups act isometrically on
each $g_{r_h}$, and $\sigma$ normalises $H$.

\begin{theorem}\label{thm:main}
Let $g_{B_1}$ be a B\"ohm metric on $B_1=S^3\times S^2$ far enough in B\"ohm's sequence. There are a radius $r_->0$ (defined in
section \ref{sec:radial}: it is the radius at which a tensor mode of the black holes \eqref{bh} becomes a zero mode of
$\Delta_L+2n$) and a family of Poincar\'e--Einstein metrics $g_v$, $0<|v|$ small, on $M=\mathbb R^2\times S^3\times S^2\times S^2$, smooth and
smoothly conformally compact, invariant under $H$ but not under $\sigma$, with $g_v\to g_{r_-}$ as $v\to0$ and conformal
infinity $[d\tau^2\oplus g_{B_1}\oplus g_{S^2}]$ (for a period $\tau_0$ close to that of $g_{r_-}$). The Killing fields of
$S^2$ which are not in the Lie algebra of $K$ do not extend to Killing fields of $g_v$.
\end{theorem}

By \eqref{pi1}, theorem \ref{thm:main} gives a counterexample to \cite[Th.~5.4]{An2}, according to which the
extension holds for every Poincar\'e--Einstein metric with $\pi_1(M,\partial M)=0$; in section \ref{sec:anderson}
we indicate the step of its proof which, in our reading, does not apply in general.

To our knowledge, in all previously known Poincar\'e--Einstein metrics with $\pi_1(M,\partial M)=0$, the Killing fields of the conformal
infinity extend to Killing fields of the bulk: this holds for the non-degenerate ones (see e.g. \cite[Prop.~A.1]{ACD2}),
for instance for hyperbolic space \cite{GL} and, more generally, whenever the sectional curvature is non-positive
\cite[Th.~A]{Lee}, in particular for those obtained by perturbation, and the explicit or cohomogeneity one examples are invariant under the
isometries of their conformal infinity. The metrics of theorem \ref{thm:main} seem to be the first Poincar\'e--Einstein metrics
with $\pi_1(M,\partial M)=0$ which have less symmetry than their conformal infinity.

\medskip\noindent\emph{Strategy of the proof.} The proof proceeds in six steps.
\begin{enumerate}[label=\arabic*.]
\item \emph{The horizon factor} (section \ref{sec:bohm}). Far enough in B\"ohm's sequence, the lowest eigenvalue
$\lambda_1$ of $\Delta_L$ on $\hat G_1$-invariant TT tensors of $B_1$ is simple, very negative, and separated from the
next one by a large gap (theorem \ref{thm:bohm} below).
\item \emph{A charged horizon mode} (section \ref{sec:horizon}). If $h_1$ is the corresponding eigentensor and $Y$ the
zonal first spherical harmonic of $S^2$, then $Yh_1$ is a TT eigentensor of $\Delta_L$ on $B$, with eigenvalue
$\lambda_*=\lambda_1+2E<0$, which is $H$-invariant and odd under $\sigma$; on $H$-invariant tensors of $B$, the spectrum of $\Delta_L$ is bounded
below by $\lambda_1$ and the eigenspace of $\lambda_*$ is $\mathbb RYh_1$ (lemma \ref{lem:bottom}).
\item \emph{The radial problem} (section \ref{sec:radial}). The tensor modes $r^2\psi(r)Yh_1$ of $g_{r_h}$ lead to a
radial Sturm--Liouville operator, whose lowest eigenvalue $\nu(r_h)$ changes sign at a radius $r_-$, the bifurcation
radius; by comparison with the planar black brane (the limit of the black holes $g_{r_h}$ as
$r_h\to\infty$, in the variables $r/r_h$ and $r_h\tau$, whose horizon is the flat $\mathbb R^N$, see the proof of lemma
\ref{lem:cross}), no other tensor mode degenerates at $r_-$.
\item \emph{The other sectors} (section \ref{sec:scalar}). The $H$-invariant kernel of $\Delta_L+2n$ at $r_-$ contains
no modes of vector or scalar type: vector modes are excluded by symmetry and gauge arguments, modes built on
constant functions by the linearised Birkhoff theorem (they would force $r_-^2=\frac{N-1}{N+1}$, which is excluded
because $r_-$ is large), and the other scalar modes by theorem \ref{thm:scalar} below.
Hence this kernel is spanned by the $\sigma$-odd tensor $k_0=r^2\psi_0Yh_1$ (section \ref{sec:kernel}).
\item \emph{Bifurcation} (section \ref{sec:bif}). A Lyapunov--Schmidt reduction of the Einstein equation in Bianchi
gauge, in the space of $H$-invariant tensors, gives an odd reduced equation, since $k_0$ is $\sigma$-odd; the change of
sign of $\nu$ at $r_-$ produces a branch of Einstein metrics $g_v$ which are not $\sigma$-invariant, and which are
smooth and smoothly conformally compact.
\item \emph{Non-extension} (section \ref{sec:bif}). Suppose that the rotation $X_0$ about the $x$-axis extends to
a Killing field $X_v$ of $g_v$. Then the identity
$\mathcal L_{X_0}k(v)=-\mathcal L_{X_v-X_0}g_v$, paired with $\hat k=\mathcal L_{X_0}k_0$, gives a contradiction at
first order in $v$.
\end{enumerate}
Section \ref{sec:anderson} contains remarks on the extension problem, on weaker hypotheses, on the Lorentzian
case and on dimension four.

Steps 1 and 4 rely on two results of independent interest. The first one concerns the B\"ohm metrics $g_i$ themselves.

\begin{theorem}\label{thm:bohm}
On $G_1$-invariant TT tensors of a B\"ohm metric $g_i$ on $S^3\times S^2$, normalised by $\Ric(g_i)=Eg_i$ for a
fixed $E>0$, every eigenvalue of $\Delta_L$ different from $0$ and $\frac54E$ is simple. Let $\lambda_1^+(g_i)<\lambda_2^+(g_i)$ be the two lowest eigenvalues of $\Delta_L$ on $\hat G_1$-invariant TT
tensors, and $\epsilon_i:=b_i(0)$ the radius of the non-collapsing sphere at the singular orbits, which tends to $0$
along B\"ohm's sequence. Then
$$
\epsilon_i^2\lambda_1^+(g_i)\to\lambda_1^\infty,\qquad \epsilon_i^2\lambda_2^+(g_i)\to\lambda_2^\infty,\qquad
\lambda_1^\infty<\lambda_2^\infty<0,
$$
where $\lambda_1^\infty,\lambda_2^\infty$ are the two lowest eigenvalues of $\Delta_L$ on $G_1$-invariant $L^2$ tensors of
B\"ohm's complete Ricci-flat metric $g_\infty$ on $\mathbb R^3\times S^2$, normalised so that its singular orbit $S^2$ has
radius $1$. In particular
$\lambda_1^+(g_i)\to-\infty$ and $\lambda_2^+(g_i)-\lambda_1^+(g_i)\to+\infty$ along B\"ohm's sequence.
\end{theorem}

Only the exceptional values $0$ and $\frac54E$ depend on the normalisation: $\epsilon_i^2\lambda$ is scale invariant,
so the limits above are not affected by a homothety. The simplicity comes from a count of regular solutions at the singular orbit. The asymptotics come from a rescaling
at the singular orbits, where the B\"ohm metrics converge to B\"ohm's complete Ricci-flat metric on
$\mathbb R^3\times S^2$ \cite{Bo2}, which is asymptotic to the cone over $S^2\times S^2$; this cone is unstable, by the
optimality of the constant in Hardy's inequality, which gives the negative eigenvalues of the model, and a curvature
bound $|\mathrm{Rm}|\le C t^{-2}$, valid for every cohomogeneity one Einstein metric of this type (lemma
\ref{lem:curv}), controls the rest of the manifold. The fact that $\Delta_L$ has negative eigenvalues on the B\"ohm
metrics was observed in \cite{GHP}, and their instability was established by Verger \cite{Ve}, using their convergence
to the sine cone \cite[Th.~3.7]{Bo}; the test tensors used here are of the type considered in \cite[\S6]{Ve} (lemma \ref{lem:quot}),
but we do not use this convergence.

The radial analysis of step 3 uses the Gibbons--Hartnoll mechanism \cite{GH} and a comparison with the planar black
brane. In step 4, vector modes are excluded by the symmetry in $\tau$ together with a gauge argument (lemma
\ref{lem:vector}); modes built on constant functions are treated by the linearised Birkhoff theorem \cite{Kl}; and for
the other scalar modes we use the master equation of Kodama and Ishibashi \cite{KI,KI2}, in the Riemannian form of
\cite{CDK}, and prove:

\begin{theorem}\label{thm:scalar}
Let $(B^7,g_B)$ be a compact Einstein manifold with $\Ric_B=6g_B$ and $r_h\ge6$. Then the static scalar-type master
equation of Kodama and Ishibashi for the black hole \eqref{bh} \textup{(}recalled in \eqref{master}\textup{)} has no
non-zero solution which is bounded at the horizon and normalisable at infinity, for every eigenvalue
$\kappa>7$ of the Laplacian $\Delta=\delta d$ of $(B,g_B)$ acting on functions. If $(B,g_B)$ is not the round sphere,
this holds for every eigenvalue $\kappa>0$.
\end{theorem}

The proof uses an explicit positive supersolution suggested by the Bessel functions of the asymptotic region, and
reduces to the positivity of an explicit polynomial. After a suitable change of variables, this polynomial is a
sum of terms whose positivity is elementary; the computations are detailed in appendix \ref{app:pot}. The same method shows that, for the planar
black brane in dimension nine, the static master equation of scalar type has no non-zero solution which is bounded at the
horizon and normalisable at infinity, for every momentum $q>0$ (appendix \ref{app:pot}, A.5). Non-degeneracy
of Riemannian Schwarzschild--anti-de Sitter metrics has been studied in \cite{CDK,Kl}, for maximally symmetric horizons
and mostly small masses; theorem \ref{thm:scalar} concerns large masses and Einstein horizons.

\section{Notation and conventions}\label{sec:notation}

\emph{Conventions.} $(\div h)_i=-\nabla^jh_{ij}$, $\delta=\div$ on symmetric $2$-tensors and on
$1$-forms, $\delta^*$ its formal adjoint, $(\delta^*\omega)_{ij}=\frac12(\nabla_i\omega_j+\nabla_j\omega_i)$, so that
$\mathcal L_Xg=2\delta^*X^\flat$, and $\mathring{\mathcal L}\omega$ is the trace-free part of $\delta^*\omega$.
$\Delta=\delta d\ge0$ on functions, $\Delta_H=d\delta+\delta d\ge0$ on $1$-forms, and the Lichnerowicz Laplacian is
$\Delta_Lh=\nabla^*\nabla h+\Ric\circ h+h\circ\Ric-2\mathring Rh$, $(\mathring Rh)_{ij}=R_{ikjl}h^{kl}$, normalised by
$\mathring Rg=\Ric$; on an Einstein manifold $\Ric=Eg$ its quadratic form is
$\langle\Delta_Lh,h\rangle_{L^2}=\int|\nabla h|^2+2E|h|^2-2\langle\mathring Rh,h\rangle$. A symmetric $2$-tensor $h$ is TT
(transverse traceless) if $\mathrm{tr}\,h=0$ and $\delta h=0$. $\ker_{L^2}(\Delta_L+2n)$ is the space of $L^2$ tensors
$k$ with $(\Delta_L+2n)k=0$, and $\beta(k)=\delta k+\frac12d\,\mathrm{tr}\,k$ is the Bianchi operator. A prime denotes
$d/dr$, a dot $d/dt$, and a subscript $u$ denotes $\partial/\partial u$ (sections \ref{sec:scalar} and
\ref{app:pot}): thus $\mathcal H_u=d\mathcal H/du$.

\emph{Spherical harmonics.} Let $\xi=(\xi^1,\dots,\xi^m)$ be the Euclidean coordinates of $\mathbb R^m$. A
\emph{spherical harmonic of degree $\ell$} on the unit sphere $S^{m-1}\subset\mathbb R^m$ is
the restriction to $S^{m-1}$ of a homogeneous harmonic polynomial of degree $\ell$ in $\xi$; equivalently
$\Delta Y=\ell(\ell+m-2)Y$ on $S^{m-1}$, the polynomial being then $|\xi|^\ell Y(\xi/|\xi|)$. On $S^2$ the eigenvalue is
$\ell(\ell+1)$ and the harmonics of degree $\ell$ form a space of dimension $2\ell+1$: the constants for $\ell=0$, the
restrictions of the linear forms for $\ell=1$, those of the trace-free quadratic forms for $\ell=2$. For
$u=f(|\xi|)Y$ with $Y$ of degree $\ell$ on $S^2$,
$$
\Delta_{\mathbb R^3}u=-\Big(\ddot f+\frac2t\dot f-\frac{\ell(\ell+1)}{t^2}f\Big)Y,\qquad t=|\xi|,
$$
the harmonic solutions being $f=t^\ell$ and $f=t^{-\ell-1}$. A tensor on $\mathbb R^3$ is said to be built on harmonics of
degree $\ell$ when its components in the coordinates $\xi^i$ are of the above form; thus $\delta_{ij}$ has degree $0$ and
$\xi_i\xi_j-\frac13|\xi|^2\delta_{ij}$ degree $2$. This is used at the singular orbit in the proof of lemma
\ref{lem:simple}, where the two exponents above are the indicial roots of the radial system. The eigenvalue depends on
the scale: on a sphere of curvature $K$, that is with $\Ric=Kg$, one has $\Delta Y_\ell=\ell(\ell+1)KY_\ell$, so that on
the factor $S^2$ of the horizon, normalised by $\Ric=Eg$, the first harmonics satisfy $\Delta^{S^2}Y=2EY$ (section
\ref{sec:horizon}). From section \ref{sec:radial} on, the harmonics of $(B,g_B)$ are the scalar, vector and tensor
harmonics in the sense of \cite[\S2]{KI2}, which are not of this elementary type.

\emph{Main notation.} The following objects are used throughout the paper.

\begin{longtable}{p{0.27\textwidth}p{0.54\textwidth}p{0.08\textwidth}}
$N=7$, $n=N+1=8$ & dimensions of the horizon $B$ and of the conformal infinity $\partial M$ (so $\dim M=n+1=9$) & \S\ref{sec:intro}\\
$E=N-1=6$ & Einstein constant of $B$, hence of its factors $B_1$ and $S^2$ (in section \ref{sec:bohm}, of the metric considered) & \S\ref{sec:intro}\\
$B=B_1\times S^2$, $g_B=g_{B_1}\oplus g_{S^2}$ & horizon, $B_1=S^3\times S^2$, $g_{B_1}$ a B\"ohm metric & \S\ref{sec:intro}\\
$g_i$, $a_i$, $b_i$, $T_i$ & B\"ohm metrics $dt^2+a_i^2g^{S^2}\oplus b_i^2g^{S^2}$ on $B_1$, where $t\in[0,T_i]$ and $T_i$ is the distance between the two singular orbits & \S\ref{sec:intro}\\
$G_1$, $\rho$, $\hat G_1$ & $SO(3)\times SO(3)$, the symmetry $t\mapsto T-t$, $G_1\times\langle\rho\rangle$ & \S\ref{sec:intro}\\
$\Sigma$, $K$, $\sigma$ & $SO(3)$ acting on $S^2$, rotations about the $z$-axis, rotation of angle $\pi$ about the $x$-axis & \S\ref{sec:intro}\\
$H$ & the group $K\times\hat G_1\times O(2)_\tau$ & \S\ref{sec:intro}\\
$g_{r_h}$, $f$, $r_h$, $\tau$, $\tau_0$ & black holes \eqref{bh}, their potential, horizon radius, Euclidean time and its period & \S\ref{sec:intro}\\
$X_0$, $X$ & boundary Killing field (in section \ref{sec:bif}, the rotation about the $x$-axis), its extension solving $(\nabla^*\nabla+n)X^\flat=0$ & \S\ref{sec:intro}\\
$\lambda_c=\frac{(N-1)(9-N)}4$ & threshold of \cite{GH}: a horizon eigenvalue below it produces a zero mode of $\Delta_L+2n$ & \S\ref{sec:intro}\\
$\lambda_1^+$, $\lambda_2^+$ & two lowest eigenvalues of $\Delta_L$ on the $G_1$-invariant, $\rho$-even TT tensors of $g_i$, that is on the $\hat G_1$-invariant ones (whence the $+$) & \S\ref{sec:intro}\\
$\epsilon_i$, $g_\infty$, $\lambda_j^\infty$ & $b_i(0)$, B\"ohm's Ricci-flat metric on $\mathbb R^3\times S^2$ with $b_\infty(0)=1$, its two lowest eigenvalues of $\Delta_L$ on $G_1$-invariant $L^2$ tensors & \S\ref{sec:intro}\\
$\mathsf M_j$ & min-max values of the quadratic form of $\Delta_L$ on invariant tensors & \S\ref{sec:bohm}\\
$\lambda_1$, $\lambda_2$, $h_1$ & $\lambda_1^+(g_i)$, $\lambda_2^+(g_i)$ for $\Ric=6g$, and an eigentensor for $\lambda_1$ & \S\ref{sec:horizon}\\
$Y$, $\lambda_*$ & zonal first spherical harmonic of $S^2$, $\lambda_*=\lambda_1+2E$ & \S\ref{sec:horizon}\\
$u$ & $1/r$, defining function of the conformal infinity & \S\ref{sec:radial}\\
$A_\lambda$, $\mu_j(r_h,\lambda)$ & radial operator of the tensor modes and its min-max values (its eigenvalues below the essential spectrum) & \S\ref{sec:radial}\\
$r_-$, $\nu$, $\psi_{0,r_h}$, $k_{r_h}$, $k_0$ & bifurcation radius, $\nu(r_h)=\mu_0(r_h,\lambda_*)$, ground state of $A_{\lambda_*}$ ($\psi_0:=\psi_{0,r_-}$), $k_{r_h}=r^2\psi_{0,r_h}Yh_1$, $k_0=k_{r_-}$ & \S\ref{sec:radial}\\
$\kappa$, $S$, $m$ & eigenvalue of $\Delta_B$, eigenfunction, $m=\kappa-N$ & \S\ref{sec:scalar}\\
$\Phi$, $V_S$, $\mathsf x$, $\mathcal H$, $U$ & master function, its potential $V_S=fU/(16r^2\mathcal H^2)$, $\mathsf x=r_h^{N-1}(1+r_h^2)r^{1-N}$, $\mathcal H=m+\frac{N(N+1)}2\mathsf x$, and $U$ the polynomial of \cite[(3.7)]{KI} & \S\ref{sec:scalar}\\
$r_*$, $\hat F$ & tortoise coordinate, $\hat F=fu^2$, so that $d/dr_*=-\hat F\,d/du$ & \S\ref{sec:scalar}\\
$\Psi_\varpi$, $w$, $\tilde V$ & supersolution, $w=\sqrt{9+(m+\varpi)u^2}$, modified potential & \S\ref{sec:scalar}\\
$\vartheta$, $\mathsf y$, $\mathsf s$ & $w-3$, $m+\varpi=m+5$, $\mathsf xu^2$ & \S\ref{sec:scalar}\\
$v$, $r(v)$, $k(v)$, $g_v$ & bifurcation parameter, radius, perturbation, bifurcating metrics & \S\ref{sec:bif}\\
\end{longtable}

\section{The B\"ohm metrics}\label{sec:bohm}

In this section $E$ denotes the Einstein constant of the metric under consideration, the B\"ohm metrics being
normalised by $\Ric=Eg$ for a fixed $E>0$. We first collect some
tools: the quadratic form of $\Delta_L$ on invariant test tensors, a weighted Hardy inequality, and the fact that the
negative spectrum only sees TT tensors.

\emph{Invariant tensors.} Let $g=dt^2+a^2g^{S^2}\oplus b^2g^{S^2}$ be a $G_1$-invariant metric on
$(t_0,t_1)\times S^2\times S^2$. Since $S^2=SO(3)/SO(2)$ carries no invariant $1$-form, and there is no invariant
mixed tensor between the two spheres, a $G_1$-invariant symmetric $2$-tensor reads
$$
k=\phi_0\,dt^2+\phi_1a^2g^{S^2}\oplus\phi_2b^2g^{S^2},
$$
where $\phi_0,\phi_1,\phi_2$ are functions of $t$. Its trace is $\phi_0+2\phi_1+2\phi_2$ and, with $\alpha_1:=\dot a/a$, $\alpha_2:=\dot b/b$, a direct
computation (using $\nabla_{\partial_t}\partial_t=0$, $\nabla_e\partial_t=\alpha_1e$ for $e$ tangent to the first sphere and
$\alpha_2e$ for $e$ tangent to the second one) gives
$$
\delta k=-\big(\dot\phi_0+2\alpha_1(\phi_0-\phi_1)+2\alpha_2(\phi_0-\phi_2)\big)dt .
$$

The test tensors are those with $\phi_0=0$ and $\phi_1=-\phi_2=\eta$, for a function $\eta$ of $t$; we denote them by
$\Upsilon_{a,b}(\eta):=\eta\,a^2g^{S^2}\oplus(-\eta\,b^2g^{S^2})$. Their Rayleigh quotient is an explicit one-dimensional integral.

\begin{lemma}\label{lem:quot}
Let $g=dt^2+a^2g^{S^2}\oplus b^2g^{S^2}$ be an Einstein metric, $\Ric=Eg$, on $(t_0,t_1)\times S^2\times S^2$, let
$\eta\in C^\infty_c((t_0,t_1))$ and $\Upsilon_{a,b}(\eta):=\eta a^2g^{S^2}\oplus(-\eta b^2g^{S^2})$. Then $\Upsilon_{a,b}(\eta)$ is
trace-free, it is divergence-free if $\alpha_1=\alpha_2$ (in particular if $a=b$), and, with $\alpha_1,\alpha_2$ as above,
$$
\begin{gathered}
\|\Upsilon_{a,b}(\eta)\|^2=4(4\pi)^2\!\int\eta^2a^2b^2dt,\\
\langle\Delta_L\Upsilon_{a,b}(\eta),\Upsilon_{a,b}(\eta)\rangle\\=4(4\pi)^2\!\int\Big(\dot\eta^2+\big(2(\alpha_1-\alpha_2)^2+2E-a^{-2}-b^{-2}\big)\eta^2\Big)a^2b^2dt .
\end{gathered}
$$
\end{lemma}
\begin{proof}
The first two statements follow from the formulas above ($\phi_0=0$, $\phi_1=-\phi_2=\eta$, and $\alpha_1=\alpha_2$ if $a=b$). Let $e_0=\partial_t$,
$(e_1,e_2)$ and $(e_3,e_4)$ be orthonormal frames of the two spheres (for $a^2g^{S^2}$ and $b^2g^{S^2}$), and
$k=\Upsilon_{a,b}(\eta)$,
so that $k=\mathrm{diag}(0,\eta,\eta,-\eta,-\eta)$ and $|k|^2=4\eta^2$.
The frame is parallel along $\partial_t$, and the
only non-zero components of $\nabla k$ are $(\nabla_0k)_{jj}=\pm\dot\eta$ ($j\ge1$), $(\nabla_jk)_{0j}=(\nabla_jk)_{j0}=-\alpha_1\eta$
($j=1,2$) and $=\alpha_2\eta$ ($j=3,4$), so that $|\nabla k|^2=4\dot\eta^2+4(\alpha_1^2+\alpha_2^2)\eta^2$. Since $k$ is diagonal and $k_{00}=0$,
$\langle\mathring Rk,k\rangle=\sum_{i\ne j}\sec_{ij}k_{ii}k_{jj}$, where the sectional curvatures are
$\sec_{12}=a^{-2}-\alpha_1^2$, $\sec_{34}=b^{-2}-\alpha_2^2$ and $\sec_{jl}=-\alpha_1\alpha_2$ ($j=1,2$, $l=3,4$); hence
$\langle\mathring Rk,k\rangle=2\eta^2(a^{-2}-\alpha_1^2+b^{-2}-\alpha_2^2+4\alpha_1\alpha_2)$. The volume form is $(4\pi)^2a^2b^2dt$ after integration over the
spheres, and the formula follows from the expression of the quadratic form of $\Delta_L$ given in the conventions.
\end{proof}

The negative part of the quotient comes from the optimality of the constant in the weighted Hardy inequality
$$\int_0^\infty r^p\psi'^2dr\ge\frac{(p-1)^2}4\int_0^\infty r^{p-2}\psi^2dr,\qquad \psi\in C^\infty_c((0,\infty)),\ p>1,$$
\cite[Ch.~IX]{HLP}; we only need the following elementary statement.

\begin{lemma}\label{lem:hardy}
Let $p>1$ and $c>\frac{(p-1)^2}4$. There are $R>1$ and $\psi\in C^\infty_c((1,R))$, $\psi\ne0$, such that
$\int_0^\infty r^p\psi'^2dr<c\int_0^\infty r^{p-2}\psi^2dr$. For $\varsigma>0$, $\psi(\cdot/\varsigma)$ is supported in $(\varsigma,\varsigma R)$, and
both integrals are multiplied by $\varsigma^{p-1}$.
\end{lemma}
\begin{proof}
Put $a:=\frac{p-1}2$, so that $c>a^2$ is the constant of the statement, and fix $\varphi_1\in C^\infty_c((0,1))$,
$\varphi_1\ne0$. For $R>1$ put $\ell:=\log R$, $\varphi(\theta):=\varphi_1(\theta/\ell)\in C^\infty_c((0,\ell))$ and
$\psi(r):=r^{-a}\varphi(\log r)$, so that $\psi\in C^\infty_c((1,R))$. Then
$\psi'(r)=r^{-a-1}\big(\frac{d\varphi}{d\theta}-a\varphi\big)(\log r)$, and the change of variable $r=e^\theta$, $dr=r\,d\theta$,
gives, since $r^p\cdot r^{-2a-2}\cdot r=1$ and $r^{p-2}\cdot r^{-2a}\cdot r=1$,
$$
\int_0^\infty r^p\psi'^2dr=\int_0^\ell\Big(\frac{d\varphi}{d\theta}-a\varphi\Big)^2d\theta,\qquad
\int_0^\infty r^{p-2}\psi^2dr=\int_0^\ell\varphi^2d\theta .
$$
Expanding the square and using $\int_0^\ell\varphi\frac{d\varphi}{d\theta}d\theta=\frac12\big[\varphi^2\big]_0^\ell=0$, we get
$$
\frac{\int r^p\psi'^2dr}{\int r^{p-2}\psi^2dr}=a^2+\frac{\int_0^\ell|\frac{d\varphi}{d\theta}|^2d\theta}{\int_0^\ell\varphi^2d\theta}
=a^2+\frac1{\ell^2}\,\frac{\int_0^1|\frac{d\varphi_1}{d\theta_1}|^2d\theta_1}{\int_0^1\varphi_1^2d\theta_1},
$$
since $\int_0^\ell|\frac{d\varphi}{d\theta}|^2d\theta=\ell^{-1}\int_0^1|\frac{d\varphi_1}{d\theta_1}|^2d\theta_1$ and $\int_0^\ell\varphi^2d\theta=\ell\int_0^1\varphi_1^2d\theta_1$ ($\theta_1=\theta/\ell$).
The last term tends to $0$ as $R\to\infty$; choosing $R$ so large that it is $<c-a^2$ gives the inequality of the
statement. Finally, the change of variable $r=\varsigma r_1$ gives
$$
\int r^p\Big(\frac d{dr}\psi(r/\varsigma)\Big)^2dr=\varsigma^{p-1}\int r_1^p\psi'(r_1)^2dr_1,
$$
$$
\int r^{p-2}\psi(r/\varsigma)^2dr=\varsigma^{p-1}\int r_1^{p-2}\psi(r_1)^2dr_1 .
$$
\end{proof}

Finally, the negative spectrum of $\Delta_L$ on an Einstein manifold with $E\ge0$ only sees TT tensors, so that test
tensors need not be TT.

\begin{lemma}\label{lem:negTT}
Let $(\mathcal M,g)$ be an Einstein manifold, $\Ric=Eg$, $E\ge0$, compact or complete with bounded curvature, and $\Gamma$ a
group of isometries. If $k$ is a $\Gamma$-invariant $L^2$ eigentensor of $\Delta_L$ with eigenvalue $\lambda<0$, then $k$
is TT. Consequently the negative eigenvalues of $\Delta_L$ on $\Gamma$-invariant symmetric $2$-tensors coincide, with
their eigenspaces, with those on $\Gamma$-invariant TT tensors, and the $j$-th negative eigenvalue on $\Gamma$-invariant TT
tensors is the $j$-th min-max value of the quadratic form $q_g(k)=\langle\Delta_Lk,k\rangle$ on all $\Gamma$-invariant
tensors of $H^1$, as long as this value is negative and below the essential spectrum.
\end{lemma}
\begin{proof}
On an Einstein manifold, $\mathrm{tr}\,\Delta_L=\Delta\,\mathrm{tr}$ and $\delta\Delta_L=\Delta_H\delta$ \cite{Li61}. Hence
$\Delta\,\mathrm{tr}\,k=\lambda\,\mathrm{tr}\,k$ and $\Delta_H\delta k=\lambda\,\delta k$ (in the complete case $k\in H^2$ by
elliptic regularity, which justifies the integrations by parts); since $\Delta\ge0$ and $\Delta_H=\nabla^*\nabla+\Ric\ge0$,
$\mathrm{tr}\,k=0$ and $\delta k=0$. The rest follows from the min-max principle for the self-adjoint operator $\Delta_L$
on $\Gamma$-invariant $L^2$ tensors. Note finally the behaviour under scaling: $\Delta_L$ of $c^2g$ is $c^{-2}$ times
$\Delta_L$ of $g$, so that the eigenvalues and the min-max values get multiplied by $c^{-2}$. This is what will relate
the B\"ohm metrics to their rescaled limits in proposition \ref{prop:a}.
\end{proof}

Next we show that these eigenvalues are simple. This is a local statement at a singular orbit: we count the
solutions of the eigenvalue equation which are smooth at the orbit, and compare them with the solutions which are
not TT.

Let $g=dt^2+a(t)^2g^{S^2}\oplus b(t)^2g^{S^2}$ be a $G_1$-invariant Einstein metric, $\Ric=Eg$, $E\ge0$, defined on a
neighbourhood of a singular orbit $\{t=0\}$ where the first $S^2$ collapses. This neighbourhood is
$\mathbb R^3\times S^2$, $t$ being the distance to the origin of $\mathbb R^3$, and the smoothness of $g$ at the orbit
means that $a$ extends to a smooth odd function with $\dot a(0)=1$ and $b$ to a smooth even positive function
\cite[\S1]{EW}. For $\lambda\in\mathbb R$, let $\mathcal S_\lambda$ be the space of $G_1$-invariant symmetric $2$-tensors $k$
(not necessarily TT), smooth on a neighbourhood of $\{t=0\}$, with $\Delta_Lk=\lambda k$ there, and $\mathcal R_\lambda$ the
space of $G_1$-invariant functions $\varphi$ (i.e. functions of $t$), smooth on a neighbourhood of $\{t=0\}$, with
$\Delta\varphi=\lambda\varphi$ there.

\begin{lemma}\label{lem:simple}
Let $\lambda\notin\{0,\frac54E\}$. The space of $G_1$-invariant TT tensors, smooth near the singular orbit, satisfying
$\Delta_Lk=\lambda k$ there, is at most one-dimensional. Consequently every eigenvalue $\lambda\notin\{0,\frac54E\}$ of
$\Delta_L$ on $G_1$-invariant TT tensors is simple, on $B_1$ (in each
$\rho$-parity sector, or without $\rho$) as well as for the $L^2$ eigenvalues of a complete metric of this type.
\end{lemma}
\begin{proof}
The proof has five steps. Step 1 writes the eigenvalue equation near the orbit as a second order system of ordinary
differential equations with a regular singular point, and computes its indicial roots. Step 2 deduces that the space
$\mathcal S_\lambda$ of solutions which are smooth at the orbit, TT or not, has dimension at most $3$. Step 3 produces
a non-zero element of $\mathcal R_\lambda$, and step 4 uses it to exhibit two linearly independent non-TT elements of
$\mathcal S_\lambda$, which leaves at most one dimension for the TT ones. Step 5 passes from this local statement to
the eigenspaces.

\smallskip
\emph{Step 1: the radial system.} Write an invariant tensor as $k=\phi_0\,dt^2+\phi_1a^2g^{S^2}\oplus\phi_2b^2g^{S^2}$ as
above, and replace $(\phi_0,\phi_1)$ by
$$
\phi_{\rm tr}:=\tfrac13(\phi_0+2\phi_1),\qquad \phi_{\rm tf}:=\phi_0-\phi_1,
$$
so that
$$
\phi_0\,dt^2+\phi_1a^2g^{S^2}=\phi_{\rm tr}\,(dt^2+a^2g^{S^2})+\phi_{\rm tf}\,\big(dt^2-\tfrac13(dt^2+a^2g^{S^2})\big).
$$
These are functions of $t$. Their meaning is seen on the model: if $a(t)=t$, then $dt^2+a^2g^{S^2}$ is the Euclidean
metric $\delta_{ij}d\xi^id\xi^j$ of $\mathbb R^3$, with $t=|\xi|$, and $dt^2=\frac{\xi_i\xi_j}{|\xi|^2}d\xi^id\xi^j$. Thus
$\phi_{\rm tr}\,\delta$ is a function times the metric, built on the spherical harmonics of degree $\ell=0$ of the
collapsing sphere, while $\phi_{\rm tf}\,(\frac{\xi_i\xi_j}{|\xi|^2}-\frac13\delta_{ij})$ is trace-free, built on those of
degree $\ell=2$; and $\phi_2$ multiplies the metric of the second sphere ($\ell=0$ again).

Put $\Phi=(\phi_{\rm tr},\phi_{\rm tf},\phi_2)$. By the formulas for the connection and the curvature of $g$ (see the proof
of lemma \ref{lem:quot}), $\Delta_Lk=\lambda k$ is a linear system of second order ordinary differential equations for
$\Phi$, whose coefficients are polynomials in $\dot a/a$, $\dot b/b$, $\ddot a/a$, $\ddot b/b$, $a^{-2}$, $b^{-2}$ and
$\lambda$. Multiplying by $-t^2$, it reads
\begin{equation}\label{radsys}
t^2\ddot\Phi+t\,P(t)\dot\Phi+Q(t)\Phi=0,
\end{equation}
where $P$ and $Q$ are $3\times3$ matrices depending smoothly on $t\in[0,\varepsilon)$: indeed, since $a$ is odd with
$\dot a(0)=1$ and $b$ is even and positive, $t\dot a/a$, $t^2/a^2$, $t^2\ddot a/a$, $t\dot b/b$, $t^2\ddot b/b$ and $t^2/b^2$
are smooth at $t=0$. Moreover $t\dot a/a=1+O(t^2)$ and $t^2/a^2=1+O(t^2)$, while the terms involving $b$ or $\ddot a$ are
$O(t^2)$ after multiplication by $t$ or $t^2$; hence $P(0)$ and $Q(0)$ are those of the product model
$\mathbb R^3\times S^2$ (flat $\mathbb R^3$ times a round sphere), on which $\Delta_L$ acts on the three components as the Euclidean Laplacian on functions of
$|\xi|$ times spherical harmonics of degree $\ell=0$, $2$ and $0$ (the second factor being parallel). Since this
Laplacian is $-(\ddot f+\frac2t\dot f-\frac{\ell(\ell+1)}{t^2}f)$ on $f(t)$ times a harmonic of degree $\ell$,
$$
P(0)=2\,I,\qquad Q(0)=-\,\mathrm{diag}(0,6,0).
$$
The point $t=0$ is a regular singular point (a singularity of the first kind) of \eqref{radsys} \cite[Ch.~4]{CL}. Its
indicial equation is $\det\big(s(s-1)I+sP(0)+Q(0)\big)=0$, i.e. $(s-\ell)(s+\ell+1)=0$ for $\ell=0,2,0$: the indicial
roots are $\{0,-1\}$, $\{2,-3\}$ and $\{0,-1\}$, namely the exponents of $|\xi|^\ell$ and $|\xi|^{-\ell-1}$.

\smallskip
\emph{Step 2: $\dim\mathcal S_\lambda\le3$.} Put $Z=(\Phi,t\dot\Phi)\in\mathbb R^6$. Then \eqref{radsys} becomes the first
order system $t\dot Z=A(t)Z$ with
$$
A(t)=\begin{pmatrix}0&I\\-Q(t)&I-P(t)\end{pmatrix},
$$
whose eigenvalues at $t=0$ are the six indicial roots. In the variable $\theta=-\log t$, it reads
$\frac{dZ}{d\theta}=-A(e^{-\theta})Z$, with $-A(e^{-\theta})=-A(0)+O(e^{-\theta})$ as $\theta\to\infty$. The eigenvalues of
$-A(0)$ are $0,0,-2$ and $1,1,3$; they are separated by the line $\mathrm{Re}=\frac12$, and the perturbation is
integrable. By the roughness of exponential dichotomies (see \cite[Ch.~4]{Cop}, or the theory of asymptotic integration
of linear systems \cite[Ch.~X]{Har}), the solutions with $|Z|=O(e^{\theta/2})$ as $\theta\to\infty$ form a space of
dimension $3$, the number of eigenvalues of $-A(0)$ with real part $<\frac12$. If $k\in\mathcal S_\lambda$, $\Phi$ is smooth
on $[0,\varepsilon)$, so $Z$ is bounded, and $k\mapsto Z$ is injective. Hence $\dim\mathcal S_\lambda\le3$. (By the
Frobenius method \cite[Ch.~4]{CL}, one may also check that $\dim\mathcal S_\lambda=3$, one solution for each
non-negative indicial root, but we do not need it.)

\smallskip
\emph{Step 3: $\dim\mathcal R_\lambda\ge1$.} For $\varepsilon$ small, $\lambda$ is not a Dirichlet eigenvalue of $\Delta$ on the
tube $U_\varepsilon=\{t<\varepsilon\}$, since the lowest Dirichlet eigenvalue of $U_\varepsilon$ tends to $+\infty$ as
$\varepsilon\to0$. The Dirichlet problem $\Delta\varphi=\lambda\varphi$ in $U_\varepsilon$, $\varphi=1$ on $\partial U_\varepsilon$, has
therefore a unique solution: writing $\varphi=1+\varphi_1$, it reads $(\Delta-\lambda)\varphi_1=\lambda$ with
$\varphi_1=0$ on the boundary, and $\lambda$ is not a Dirichlet eigenvalue. This solution is smooth in $U_\varepsilon$
by elliptic regularity, $G_1$-invariant by uniqueness since the domain and the data are, and non-zero since it equals
$1$ on the boundary: $0\ne\varphi\in\mathcal R_\lambda$. Only the existence of such a $\varphi$ is used below (in fact
$\dim\mathcal R_\lambda=1$: the scalar equation has the indicial roots $0$ and $-1$, and step 2 applies to it as
well). Note also that $\mathcal R_\lambda$ contains no non-zero constant, since $\lambda\ne0$.

\smallskip
\emph{Step 4: the non-TT solutions.} On an Einstein manifold, $\mathrm{tr}\,\Delta_L=\Delta\,\mathrm{tr}$,
$\delta\Delta_L=\Delta_H\delta$, $\Delta_Hd=d\Delta$, $\Delta_L(\varphi g)=(\Delta\varphi)g$ and $\Delta_L\delta^*=\delta^*\Delta_H$,
hence $\Delta_L\nabla^2\varphi=\nabla^2\Delta\varphi$ since $\nabla^2\varphi=\delta^*d\varphi$ \cite{Li61}. Let $k\in\mathcal S_\lambda$.
Then $\mathrm{tr}\,k\in\mathcal R_\lambda$. The $1$-form $\delta k$ is $G_1$-invariant; the invariant $1$-forms being the $\varsigma(t)\,dt$, it is closed, hence
exact on the simply connected neighbourhood $\mathbb R^3\times S^2$: $\delta k=d\tilde\varphi$ with
$\tilde\varphi$ a smooth function of $t$ alone, hence $G_1$-invariant. Since $d(\Delta\tilde\varphi-\lambda\tilde\varphi)=\Delta_H\delta k-\lambda\delta k=0$, the function
$\Delta\tilde\varphi-\lambda\tilde\varphi$ is a constant $c$, and replacing $\tilde\varphi$ by $\tilde\varphi+c/\lambda$ we get
$\tilde\varphi\in\mathcal R_\lambda$; this normalisation determines $\tilde\varphi$, since $\mathcal R_\lambda$ contains no non-zero
constant. This defines a linear map
$$
\Pi:\mathcal S_\lambda\to\mathcal R_\lambda\times\mathcal R_\lambda,\qquad k\mapsto(\mathrm{tr}\,k,\tilde\varphi),
$$
and $\Pi(k)=0$ if and only if $\mathrm{tr}\,k=0$ and $\delta k=0$: the kernel of $\Pi$ is the space of TT solutions.

Let $0\ne\varphi\in\mathcal R_\lambda$ (step 3). Then $\varphi g$ and $\nabla^2\varphi$ belong to $\mathcal S_\lambda$. Since
$\dim B_1=5$, $\mathrm{tr}(\varphi g)=5\varphi$ and $\mathrm{tr}\,\nabla^2\varphi=-\Delta\varphi=-\lambda\varphi$. With our convention
$\delta=\div$, $\delta(\varphi g)=-d\varphi$, and the Weitzenb\"ock formula
$2\delta\delta^*\omega=\nabla^*\nabla\omega+d\delta\omega-\Ric(\omega)$ together with $\nabla^*\nabla d\varphi=\Delta_Hd\varphi-\Ric(d\varphi)$
gives $2\delta\nabla^2\varphi=2\delta\delta^*d\varphi=2(\lambda-E)d\varphi$. Hence
$$
\Pi(\varphi g)=(5\varphi,-\varphi),\qquad \Pi(\nabla^2\varphi)=(-\lambda\varphi,(\lambda-E)\varphi).
$$
These two vectors are linearly independent, since $4\lambda-5E\ne0$. Hence $\Pi$ has rank at least $2$,
and its kernel has dimension at most $\dim\mathcal S_\lambda-2\le1$ by step 2.

\smallskip
\emph{Step 5: conclusion.} An invariant eigentensor of $\Delta_L$ on $B_1$, or an $L^2$ eigentensor on a complete metric
of this type, is smooth at the singular orbit, hence its restriction near the orbit is a TT element of
$\mathcal S_\lambda$. It is determined by this restriction, since on $\{t>0\}$ the eigenvalue equation is a linear system
of ordinary differential equations with smooth coefficients, for which the Cauchy problem has a unique solution.
Hence the eigenspaces have dimension at most $1$.
\end{proof}

To compare the B\"ohm metrics with their limits, near and away from the singular orbits, we need a curvature bound
which is uniform along the sequence. It holds for every cohomogeneity one Einstein metric of this type, and follows
from the Einstein equations alone.

\begin{lemma}\label{lem:curv}
Let $g=dt^2+a^2g^{S^2}\oplus b^2g^{S^2}$ on $(0,T)$ be any metric with $\Ric=Eg$, $E>0$, smooth at $t=0$ where $a(0)=0$, $\dot a(0)=1$,
$b(0)>0$, and let $L=2\dot a/a+2\dot b/b$ be the mean curvature of the orbits. On any interval $(0,t_1)$ where $L>0$,
$$
L\le\frac4t\qquad\text{and}\qquad|\mathrm{Rm}_g|\le C\,(L^2+E)\le\frac{C_1}{t^2},
$$
with universal constants $C,C_1$.
\end{lemma}
\begin{proof}
Let $\alpha_1=\dot a/a$, $\alpha_2=\dot b/b$ as above. The Einstein equations read
$$
\frac{\ddot a}a=-E-\alpha_1^2-2\alpha_1\alpha_2+\frac1{a^2},\qquad\frac{\ddot b}b=-E-\alpha_2^2-2\alpha_1\alpha_2+\frac1{b^2},\qquad\frac{\ddot a}a+\frac{\ddot b}b=-\frac E2,
$$
the first two being the components of $\Ric=Eg$ along the two spheres and the third one its component along
$\partial_t$. Adding the first two and using the third one gives the constraint
$2\alpha_1^2+2\alpha_2^2+8\alpha_1\alpha_2=\frac2{a^2}+\frac2{b^2}-3E$. Writing $\alpha_1=\frac L4+\upsilon$, $\alpha_2=\frac L4-\upsilon$, we have
$2\alpha_1^2+2\alpha_2^2=\frac{L^2}4+4\upsilon^2$ and $8\alpha_1\alpha_2=\frac{L^2}2-8\upsilon^2$, so that the constraint becomes
$$
\frac2{a^2}+\frac2{b^2}+4\upsilon^2=\frac34L^2+3E .
$$
The three terms of the left-hand side are non-negative, hence each of them is bounded by the right-hand side:
$$
\frac1{a^2}\le\frac38L^2+\frac32E,\qquad\frac1{b^2}\le\frac38L^2+\frac32E,\qquad\upsilon^2\le\frac3{16}L^2+\frac34E,
$$
and therefore $\alpha_1^2,\alpha_2^2\le2\big(\frac{L^2}{16}+\upsilon^2\big)\le\frac12L^2+\frac32E$ and $|\alpha_1\alpha_2|\le\frac12L^2+\frac32E$. By the Einstein
equations above, $|\ddot a/a|$ and $|\ddot b/b|$ are then also bounded by $C(L^2+E)$. Now the curvature operator of $g$
is diagonal in the adapted frame of the proof of lemma \ref{lem:quot}, where the sectional curvatures of the planes
tangent to the orbits are computed (see also \cite[\S1]{EW} for the curvature of these metrics), with eigenvalues $-\ddot a/a$, $-\ddot b/b$ (planes containing $\partial_t$),
$(1-\dot a^2)/a^2=a^{-2}-\alpha_1^2$, $b^{-2}-\alpha_2^2$ (planes tangent to one sphere) and $-\alpha_1\alpha_2$ (mixed planes); all of them
are bounded by $C(L^2+E)$ by the above, which gives $|\mathrm{Rm}_g|\le C(L^2+E)$.

It remains to bound $L$. Since $\dot\alpha_1=\ddot a/a-\alpha_1^2$ and $\dot\alpha_2=\ddot b/b-\alpha_2^2$, the third Einstein equation gives the
Riccati equation
$$
\dot L=2\Big(\frac{\ddot a}a+\frac{\ddot b}b\Big)-2\alpha_1^2-2\alpha_2^2=-E-2\alpha_1^2-2\alpha_2^2\le-E-(\alpha_1+\alpha_2)^2=-E-\frac{L^2}4 .
$$
Hence, where $L>0$, $\frac{d}{dt}\big(\frac1L\big)=-\frac{\dot L}{L^2}\ge\frac14$. Since $\alpha_1=\dot a/a\sim1/t$ and $\alpha_2\to0$ as $t\to0$, we
have $L\sim2/t$ and $1/L\to0$; integrating from $0$ to $t$ gives $\frac1L\ge\frac t4$, that is $L\le\frac4t$. Finally
the Riccati equation also gives
$\frac d{dt}\arctan\frac L{2\sqrt E}=\frac{2\sqrt E\,\dot L}{4E+L^2}\le-\frac{\sqrt E}2$; since
$\arctan\frac L{2\sqrt E}\to\frac\pi2$
as $t\to0$ and $\arctan\frac L{2\sqrt E}>0$ where $L>0$, we get $t_1\le\frac\pi{\sqrt E}$, hence $E\le\frac{\pi^2}{t_1^2}\le\frac{\pi^2}{t^2}$ on $(0,t_1)$. Therefore
$L^2+E\le\frac{16+\pi^2}{t^2}$ there, and the last inequality of the statement holds with
$C_1=C(16+\pi^2)$.
\end{proof}

For the B\"ohm metrics, $L>0$ on $(0,T/2)$ ($L$ is decreasing and vanishes at the middle, by symmetry), so lemma
\ref{lem:curv} gives $|\mathrm{Rm}_{g_i}|\le C_1t^{-2}$ on $(0,T_i/2]$, uniformly in $i$.

We now rescale the B\"ohm metrics at their singular orbits. The rescaled metrics converge to a complete Ricci-flat
metric on $\mathbb R^3\times S^2$, whose Lichnerowicz Laplacian has negative eigenvalues, simple by lemma
\ref{lem:simple}; the eigenvalues of the B\"ohm metrics then follow them. We first record the continuous dependence of
the singular initial value problem on its data, which we need with respect to the Einstein constant.

\begin{lemma}\label{lem:cont}
Let $E\in\mathbb R$ and $b_0>0$. The singular initial value problem
\begin{equation}\label{ivp}
\ddot a=-Ea+\frac{1-\dot a^2}a-\frac{2\dot a\dot b}b,\qquad
\ddot b=-Eb+\frac{1-\dot b^2}b-\frac{2\dot a\dot b}a,
\end{equation}
$a(0)=0$, $\dot a(0)=1$, $b(0)=b_0$, $\dot b(0)=0$, has a unique solution, and $dt^2+a^2g^{S^2}\oplus b^2g^{S^2}$ is then
a smooth metric with $\Ric=Eg$ near the singular orbit $\{t=0\}$ \textup{(}\cite{EW}, \cite[Th.~2.3]{Bo}\textup{)}. If moreover this solution exists on $[0,R]$ with $a>0$ on $(0,R]$, then so does the solution
for every $(\tilde E,\tilde b_0)$ close to $(E,b_0)$, and the functions $\frac{a-t}{t^3}$, $\frac{\dot a-1}{t^2}$, $\frac{\dot b}t$,
$b$, $\frac{\ddot a}a$, $\frac{\ddot b}b$ converge to those of the first solution uniformly on $[0,R]$ as
$(\tilde E,\tilde b_0)\to(E,b_0)$.
\end{lemma}
\begin{proof}
Existence and uniqueness are given by the references of the statement; that the component of $\Ric=Eg$ along
$\partial_t$ follows from \eqref{ivp} for a solution which is smooth at the singular orbit is \cite[Lemma~2.1]{Bo}.
The continuous dependence on $b_0$ is part of \cite[Th.~2.3]{Bo}, and we need it with respect to the Einstein constant
as well. It follows from \cite[Th.~2.2]{Bo}: \emph{if $A$, $B$ are analytic near $y_0\in\mathbb R^j$, $B(y_0)=0$ and
$I-\frac1mDB_{y_0}$ is invertible for every integer $m\ge1$, then $\dot y=A(y)+\frac1tB(y)$ has a unique solution
with $y(0)=y_0$, and it depends continuously on $y_0\in\{B=0\}$}, uniformly on an
interval $[0,t_0]$, $t_0>0$. Indeed, let
$$
y=(y_1,\dots,y_6):=\Big(\frac{a-t}{t^3},\ \frac{\dot a-1}{t^2},\ \frac{\dot b}t,\ b,\ t,\ E\Big),
$$
where $t$ and the Einstein constant are considered as unknowns ($\dot y_5=1$, $\dot y_6=0$); for a smooth metric, $y$ is
smooth up to $t=0$. A computation gives $\dot y=A(y)+\frac1tB(y)$ with $A(y)=(0,0,0,y_5y_3,1,0)$ and
$B(y)=(y_2-3y_1,B_2(y),B_3(y),0,0,0)$,
$$
\begin{aligned}
B_2(y)&=-y_6(1+y_5^2y_1)-\frac{y_2(2+y_5^2y_2)}{1+y_5^2y_1}-\frac{2(1+y_5^2y_2)y_3}{y_4}-2y_2,\\
B_3(y)&=-y_6y_4+\frac{1-y_5^2y_3^2}{y_4}-\frac{2(1+y_5^2y_2)y_3}{1+y_5^2y_1}-y_3 ,
\end{aligned}
$$
which are analytic near the initial value $y(0)$; the latter is determined by $b_0$ and $E$ through $y_5(0)=0$ and
$B(y(0))=0$. At $y(0)$, the matrix $DB$ is upper triangular with diagonal $(-3,-4,-3,0,0,0)$, hence $I-\frac1mDB_{y(0)}$ is
invertible for every $m\ge1$, and the theorem applies: $y$ depends continuously on $(E,b_0)$, uniformly on some interval $[0,t_0]$. With the usual continuous dependence away from $t=0$, the functions $y_1,\dots,y_4$ converge uniformly on
$[0,R]$, and then so do $\ddot a/a$ and $\ddot b/b$ by \eqref{ivp}.
\end{proof}

\begin{lemma}\label{lem:rescale}
Let $\epsilon_i:=b_i(0)$ be the radius of the non-collapsing sphere at the singular orbits of $g_i$, and
$\hat g_i:=\epsilon_i^{-2}g_i$, written in $\hat t=t/\epsilon_i$. Then $\epsilon_i\to0$, $T_i/\epsilon_i\to\infty$, and for
every $R$ the metrics $\hat g_i$ converge in $C^2$ on $\{\hat t\le R\}$ to the complete Ricci-flat metric
$$
g_\infty=d\hat t^2+a_\infty^2g^{S^2}\oplus b_\infty^2g^{S^2}
$$
on $\mathbb R^3\times S^2$ with $a_\infty(0)=0$, $\dot a_\infty(0)=1$, $b_\infty(0)=1$. Moreover
$|\mathrm{Rm}_{\hat g_i}|\le C_1\hat t^{-2}$ on $(0,T_i/(2\epsilon_i)]$ and $|\mathrm{Rm}_{g_\infty}|\le C_1\hat t^{-2}$ on
$(0,\infty)$.
\end{lemma}
\begin{proof}
B\"ohm's metrics on $S^3\times S^2$ are obtained for a decreasing sequence of initial values $b_i(0)$ which converges to
$0$ \cite[\S6]{Bo}, so that $\epsilon_i\to0$. The metric $\hat g_i$ solves \eqref{ivp}
with Einstein constant $E\epsilon_i^2\to0$ and the data $a(0)=0$, $\dot a(0)=1$, $b(0)=1$, $\dot b(0)=0$, while
$g_\infty$ is by definition the solution for the Einstein constant $0$ and the same data. By \cite[Th.~6.1]{Bo2}
(in the notation of \cite{Bo2}, $r=0$, $k=2$,
$G_1/K_1=S^2$: the Ricci-flat family is then $0$-dimensional, i.e. unique up to homothety; the case $r=0$ is treated
in \cite[\S11]{Bo2}) the metric $g_\infty$ is complete, so that for every $R$ it is defined on $[0,R]$ with $a_\infty>0$ on
$(0,R]$. Lemma \ref{lem:cont} then applies: for $i$ large, $\hat g_i$ is defined on $\{\hat t\le R\}$ and the six
functions listed there converge uniformly to those of $g_\infty$. Let us check that this gives $\hat g_i\to g_\infty$ in $C^2$ on $\{\hat t\le R\}$. With $(a,b)$ the functions of
$\hat g_i$, we have $h:=\hat g_i-g_\infty=\phi_1a_\infty^2g^{S^2}\oplus\phi_2b_\infty^2g^{S^2}$, where
$\phi_1=a^2/a_\infty^2-1$ and $\phi_2=b^2/b_\infty^2-1$. Let $\alpha_1,\alpha_2$ be those of $g_\infty$. In the adapted
frame of $g_\infty$, the connection being given by $\nabla_e\partial_t=\alpha_je$ as in the proof of lemma
\ref{lem:quot}, the non-zero components of $\nabla h$ are $\dot\phi_j$ and $\alpha_j\phi_j$, and those of $\nabla^2h$ are
linear combinations of $\ddot\phi_j$, $\alpha_l\dot\phi_j$, $\dot\alpha_j\phi_j$, $\alpha_j^2\phi_j$ and
$\alpha_1\alpha_2\phi_j$. Let $\varkappa_i\to0$ be the supremum on $[0,R]$ of the differences of the six functions. Since
$a=\hat t(1+\hat t^2y_1)$ and $\dot a/a=\hat t^{-1}(1+\hat t^2y_2)/(1+\hat t^2y_1)$, we get $\phi_1=O(\hat t^2\varkappa_i)$,
$\dot\phi_1=2\frac{a^2}{a_\infty^2}\big(\frac{\dot a}a-\alpha_1\big)=O(\hat t\,\varkappa_i)$ and, using
$\frac d{d\hat t}\frac{\dot a}a=\frac{\ddot a}a-\frac{\dot a^2}{a^2}$, $\ddot\phi_1=O(\varkappa_i)$, while $\alpha_1=O(\hat t^{-1})$ and
$\dot\alpha_1=O(\hat t^{-2})$; similarly $\phi_2=O(\varkappa_i)$, $\dot\phi_2=O(\hat t\,\varkappa_i)$, $\ddot\phi_2=O(\varkappa_i)$,
with $\alpha_2=O(\hat t)$ and $\dot\alpha_2=O(1)$. Hence $|h|+|\nabla h|+|\nabla^2h|=O(\varkappa_i)$ on $\{\hat t\le R\}$. The first positive
zero of the function $a$ of $\hat g_i$ being $T_i/\epsilon_i$, this gives $T_i/\epsilon_i>R$ for $i$ large, for every $R$:
thus $T_i/\epsilon_i\to\infty$. Finally, the estimate of lemma \ref{lem:curv} being scale invariant,
$|\mathrm{Rm}_{\hat g_i}|\le C_1\hat t^{-2}$ on $(0,T_i/(2\epsilon_i)]$, and the bound for $g_\infty$ follows by letting
$i\to\infty$.
\end{proof}

The model $g_\infty$ is asymptotic to a Ricci-flat cone, and this is where its negative directions come from.

\begin{lemma}\label{lem:asympt}
With $\alpha_1=\dot a_\infty/a_\infty$, $\alpha_2=\dot b_\infty/b_\infty$ and $L=2\alpha_1+2\alpha_2$ as in lemma
\ref{lem:curv}, one has, as $\hat t\to\infty$,
$$
\begin{gathered}
\frac{a_\infty}{\hat t}\to\frac1{\sqrt3},\qquad\frac{b_\infty}{\hat t}\to\frac1{\sqrt3},\qquad
\frac{\hat t^2}{a_\infty^2}\to3,\qquad\frac{\hat t^2}{b_\infty^2}\to3,\\
\hat tL\to4,\qquad \hat t(\alpha_1-\alpha_2)\to0\,:
\end{gathered}
$$
the metric $g_\infty$ is asymptotic to the Ricci-flat cone $d\hat t^2+\frac{\hat t^2}3(g^{S^2}\oplus g^{S^2})$ over
$S^2\times S^2$.
\end{lemma}
\begin{proof}
On $\mathbb R^3\times S^2$ the cone over $(S^2\times S^2,g_E)$, $g_E$ the product Einstein metric, is the
unique cone solution, and it is stable \cite[Ex.~5.1]{Bo2}, in the sense of \cite[\S4]{Bo2}: $g_E$ is a non-degenerate
local minimum of the total scalar curvature among the invariant metrics of the same volume on the principal orbit, which
makes the cone solution a local attractor of the cohomogeneity one Einstein equations \cite[Cor.~4.2]{Bo2}. (This
dynamical notion is not the linear stability of the cone as a Ricci-flat metric, which fails here, as the negative
directions found below show.) The Ricci-flat metric $g_\infty$ is
obtained from a solution which enters one of the tubes around this cone of \cite[Prop.~4.3]{Bo2} (see
\cite[Cor.~4.4]{Bo2} and \cite[Convergence Theorem~11.1]{Bo2}). Then, by \cite[Lemma~4.5]{Bo2} (the cone
corresponds to a zero of a vector field, at which all the eigenvalues of the linearisation have a negative real part) and
\cite[Prop.~4.6]{Bo2} (whose statement records the convergence of $\bar g$, that of $d\bar g/d\tau$ being part of
its proof: the zero is asymptotically stable),
$$
\Big(\bar g,\frac{d\bar g}{d\tau}\Big)\to(\bar g_E,0)\qquad(\hat t\to\infty),
$$
where $\bar g$ is the volume normalised metric of the principal orbits and $\tau$ is the parameter of \cite[\S3]{Bo2}; here
$\bar g=\frac{a_\infty}{b_\infty}g^{S^2}\oplus\frac{b_\infty}{a_\infty}g^{S^2}$ and
$d\tau=(a_\infty b_\infty)^{-1/2}d\hat t$. With
$\alpha_1=\dot a_\infty/a_\infty$, $\alpha_2=\dot b_\infty/b_\infty$ and $L=2\alpha_1+2\alpha_2$ as in lemma
\ref{lem:curv}, we have
$\frac{d\bar g}{d\tau}=(a_\infty b_\infty)^{1/2}(\alpha_1-\alpha_2)\big(\frac{a_\infty}{b_\infty}g^{S^2}\oplus-\frac{b_\infty}{a_\infty}g^{S^2}\big)$,
so that
$$
\frac{a_\infty}{b_\infty}\to1,\qquad a_\infty(\alpha_1-\alpha_2)\to0\qquad(\hat t\to\infty).
$$
Now the constraint equation of lemma \ref{lem:curv}, for Einstein
constant $0$ instead of $E$, reads $2\alpha_1^2+2\alpha_2^2+8\alpha_1\alpha_2=\frac2{a_\infty^2}+\frac2{b_\infty^2}$, i.e., with
$\alpha_1=\frac L4+\upsilon$, $\alpha_2=\frac L4-\upsilon$,
$$
\tfrac34L^2=\frac2{a_\infty^2}+\frac2{b_\infty^2}+4\upsilon^2 .
$$
Since $a_\infty\upsilon\to0$ and $a_\infty/b_\infty\to1$, multiplying by $a_\infty^2$ we get $L\,a_\infty\to\frac4{\sqrt3}$ ($L$ is positive: it does not vanish, by this identity, and
$L\sim2/\hat t$ at $\hat t=0$; see also \cite[Prop.~3.2]{Bo2}), hence
$\alpha_1a_\infty,\alpha_2a_\infty\to\frac1{\sqrt3}$, that is $\dot a_\infty\to\frac1{\sqrt3}$ and
$\dot b_\infty=\alpha_2b_\infty\to\frac1{\sqrt3}$. In particular $a_\infty,b_\infty\to\infty$, and
$a_\infty/\hat t\to\frac1{\sqrt3}$, $b_\infty/\hat t\to\frac1{\sqrt3}$. Finally $\hat tL\to4$,
$\hat t(\alpha_1-\alpha_2)=\frac{\hat t}{a_\infty}\,a_\infty(\alpha_1-\alpha_2)\to0$ and $\hat t^2a_\infty^{-2},\hat t^2b_\infty^{-2}\to3$.

 \end{proof}

\begin{proposition}\label{prop:model}
On $G_1$-invariant $L^2$ tensors of $g_\infty$, the essential spectrum of $\Delta_L$ is contained in $[0,\infty)$, and
$\Delta_L$ has at least two negative eigenvalues. Its two lowest eigenvalues $\lambda^\infty_1<\lambda^\infty_2<0$ are
simple, and their eigentensors are TT.
\end{proposition}
\begin{proof}
By lemma \ref{lem:rescale}, $|\mathrm{Rm}_{g_\infty}|\le C_1\hat t^{-2}$; since $\Ric(g_\infty)=0$ and $|\langle\mathring Rk,k\rangle|\le
C|\mathrm{Rm}||k|^2$, we get $q_{g_\infty}(k)\ge-CC_1R^{-2}\|k\|^2$ for $k$ supported in $\{\hat t\ge R\}$, so that the
essential spectrum is contained in $[0,\infty)$ by Persson's theorem
\cite[Ch.~3]{CFKS}. 

Let $\psi\in C^\infty_c((1,R_0))$ be given by lemma \ref{lem:hardy} with $p=4$,
$c=6$. By lemmas \ref{lem:quot} ($E=0$) and \ref{lem:asympt}, the Rayleigh quotient of $\Upsilon_{a_\infty,b_\infty}(\psi(\cdot/R))$ is
$$
R^{-2}\Big(\frac{\int(r^4\psi'^2-6r^2\psi^2)dr}{\int r^4\psi^2dr}+o(1)\Big)\qquad(R\to\infty),
$$
which is negative for $R$ large. Fix $R_1$ such that this quotient is negative for $R=R_1$ and for $R=R_1R_0$. Since
$\psi$ is supported in $(1,R_0)$, the function $\psi(\cdot/R_1)$ is supported in $(R_1,R_1R_0)$ and $\psi(\cdot/(R_1R_0))$ in
$(R_1R_0,R_1R_0^2)$: the two test tensors $\Upsilon_{a_\infty,b_\infty}(\psi(\cdot/R_1))$ and
$\Upsilon_{a_\infty,b_\infty}(\psi(\cdot/(R_1R_0)))$ have disjoint supports. They are therefore orthogonal for $L^2$ and for the
quadratic form $q_{g_\infty}$, which is negative on each of them, hence negative definite on the plane they span;
by the min-max principle $\Delta_L$ has at least two negative eigenvalues
$\lambda^\infty_1\le\lambda^\infty_2<0$ on invariant $L^2$ tensors (negative modes of $\Delta_L$ on the non-compact
B\"ohm metrics were already observed in \cite{GHP}). By lemma \ref{lem:negTT} they are TT, and by lemma
\ref{lem:simple} they are simple: $\lambda^\infty_1<\lambda^\infty_2<0$. The same argument with $k$ disjointly supported
dilations shows that $\Delta_L$ has infinitely many negative eigenvalues, accumulating at $0$; we only use the two
lowest ones.

\end{proof}

We can now show that the lowest eigenvalue of the B\"ohm metrics tends to $-\infty$ and that the two lowest ones
separate along the sequence: after rescaling, they converge to the eigenvalues of the model.

\begin{proposition}\label{prop:a}
Let $\lambda^\infty_1<\lambda^\infty_2<0$ be as in proposition \ref{prop:model}. For $i$ large, the two lowest
eigenvalues $\lambda_1^+(g_i)<\lambda_2^+(g_i)$ of $\Delta_L$ on $\hat G_1$-invariant TT tensors are negative and satisfy
$$
\epsilon_i^2\lambda_j^+(g_i)\to\lambda_j^\infty,\qquad j=1,2 .
$$
In particular $\lambda_1^+(g_i)\to-\infty$ and $\lambda_2^+(g_i)-\lambda_1^+(g_i)\to+\infty$.
\end{proposition}
\begin{proof}
\emph{Step 1: reduction to min-max values.} Let $\mathsf M_1(g_i)\le\mathsf M_2(g_i)$ be the first two min-max values of the quadratic form
$q_{g_i}(k)=\int|\nabla k|^2+2E|k|^2-2\langle\mathring Rk,k\rangle$ on all $\hat G_1$-invariant tensors of $H^1$. By lemma
\ref{lem:negTT}, as soon as $\mathsf M_j(g_i)<0$ it is the $j$-th eigenvalue $\lambda_j^+(g_i)$ on $\hat G_1$-invariant TT
tensors, so that test tensors need not be TT; the same holds for $g_\infty$. We prove that $\epsilon_i^2\mathsf M_j(g_i)\to\lambda_j^\infty<0$; in particular $\mathsf M_j(g_i)<0$ for $i$ large, and
$\mathsf M_j(g_i)=\lambda_j^+(g_i)$.

\smallskip
\emph{Step 2: upper bound.} The two ends of $[0,T_i]$ are exchanged by $\rho$, and a $\rho$-even tensor is determined by
its restriction to $[0,T_i/2]$. Let $k^\infty_1,k^\infty_2$ be $L^2$-orthonormal invariant eigentensors of $g_\infty$ for
$\lambda^\infty_1,\lambda^\infty_2$, and $\chi_R$ a cut-off equal to $1$ on $\{\hat t\le R\}$, vanishing on $\{\hat t\ge2R\}$,
with $|d\chi_R|\le2/R$. Since $\chi_Rk^\infty_j\to k^\infty_j$ in $H^1$, for every $\iota>0$ and $R$ large, the $j$-th min-max value of
$q_{g_\infty}$ on $\mathrm{span}\{\chi_Rk^\infty_1,\chi_Rk^\infty_2\}$ is $\le\lambda^\infty_j+\iota$, $j=1,2$. For $i$ large,
$2R\epsilon_i<T_i/2$; transplanting $\chi_Rk^\infty_1,\chi_Rk^\infty_2$ to $\{\hat t\le2R\}$ for $\hat g_i$ and symmetrising with
$\rho$ (the two terms having disjoint supports), the $C^2$ convergence of lemma \ref{lem:rescale} gives $\rho$-even test tensors for
which the same holds with $\lambda^\infty_j+2\iota$. Since the $\rho$-even invariant tensors are the $\hat G_1$-invariant ones, and the min-max values of $q_{\hat g_i}$ are
$\epsilon_i^2\mathsf M_j(g_i)$, we get $\limsup\epsilon_i^2\mathsf M_j(g_i)\le\lambda^\infty_j$.

\smallskip
\emph{Step 3: lower bound.} Let $\chi_1^2+\chi_2^2=1$ be a $\rho$-invariant partition of unity on $B_1$, with $\chi_1=1$
where $\hat t\le R$ or $(T_i-t)/\epsilon_i\le R$, $\chi_1=0$ where both quantities are $\ge2R$, and
$|d\chi_1|^2+|d\chi_2|^2\le C/R^2$ for $\hat g_i$. On the support of $\chi_2$, $|\mathrm{Rm}_{\hat g_i}|\le C_1R^{-2}$ by
lemma \ref{lem:rescale} (and $\rho$), and the Ricci term of $q_{\hat g_i}$ is non-negative, so $q_{\hat g_i}(\chi_2k)\ge-C_2R^{-2}\|\chi_2k\|^2$.
We use the following localisation formula (the IMS formula of the theory of Schr\"odinger operators
\cite[\S3.1]{CFKS}): for a partition of unity $\chi_1^2+\chi_2^2=1$,
$$
q_{\hat g_i}(k)=q_{\hat g_i}(\chi_1k)+q_{\hat g_i}(\chi_2k)-\int(|d\chi_1|^2+|d\chi_2|^2)|k|^2 .
$$
It is elementary: $|\nabla(\chi_ak)|^2=\chi_a^2|\nabla k|^2+|d\chi_a|^2|k|^2+\frac12\langle d(\chi_a^2),d|k|^2\rangle$, the sum over
$a=1,2$ of the last terms vanishes since $\chi_1^2+\chi_2^2=1$, and the zeroth order terms of $q$ are multiplied by
$\chi_a^2$. It allows us to treat separately the neighbourhoods of the singular orbits, where $\hat g_i$ is close to
$g_\infty$, and the rest of the manifold, where the curvature is small. Thus $q_{\hat g_i}(k)+CR^{-2}\|k\|^2\ge q_{\hat g_i}(\chi_1k)+q_{\hat g_i}(\chi_2k)$, and $k\mapsto(\chi_1k,\chi_2k)$ is
an isometry of $L^2$ into the direct sum of the tensors supported in $\{\chi_1\ne0\}$ and of those supported in
$\{\chi_2\ne0\}$; the min-max values of a form on a subspace being at least those on the whole space, and the form
being $\ge-C_2R^{-2}$ on the second summand, this gives, for the $j$-th min-max value $\mathsf M_j(\hat g_i)=\epsilon_i^2\mathsf M_j(g_i)$ on $\rho$-even invariant tensors,
$$
\mathsf M_j(\hat g_i)\ge\min\big(\lambda_j^D(\hat g_i;2R),-C_2R^{-2}\big)-CR^{-2},
$$
where $\lambda_j^D(\hat g_i;2R)$ are the min-max values of $q_{\hat g_i}$ on invariant tensors supported in $\{\hat t\le2R\}$.
These depend continuously on the metric in $C^2$ on this fixed region, so $\lambda_j^D(\hat g_i;2R)\to\lambda_j^D(g_\infty;2R)
\ge\lambda^\infty_j$ (domain monotonicity and proposition \ref{prop:model}). Choosing $R$ with $C_2R^{-2}<|\lambda^\infty_2|$, we get
$\liminf\mathsf M_j(\hat g_i)\ge\lambda^\infty_j-CR^{-2}$, and $R\to\infty$ gives $\liminf\epsilon_i^2\mathsf M_j(g_i)\ge\lambda^\infty_j$.

\smallskip
\emph{Conclusion.} $\epsilon_i^2\mathsf M_j(g_i)\to\lambda^\infty_j<0$, $j=1,2$. By step 1, for $i$ large
$\mathsf M_j(g_i)=\lambda_j^+(g_i)$, hence $\epsilon_i^2\lambda_j^+(g_i)\to\lambda^\infty_j$. Since $\epsilon_i\to0$,
$\lambda_1^+(g_i)=\epsilon_i^{-2}(\lambda^\infty_1+o(1))\to-\infty$ and
$\lambda_2^+(g_i)-\lambda_1^+(g_i)=\epsilon_i^{-2}(\lambda^\infty_2-\lambda^\infty_1+o(1))\to+\infty$.
\end{proof}

\section{The spectrum of the horizon}\label{sec:horizon}

From now on $g_{B_1}=g_i$ is a B\"ohm metric on $B_1=S^3\times S^2$ normalised by $\Ric=6g_{B_1}$, so that $E=6$,
which is the normalisation of section \ref{sec:bohm} for that value of $E$. We write
$\lambda_1:=\lambda_1^+(g_i)$ and $\lambda_2:=\lambda_2^+(g_i)$, and $h_1$ for an eigentensor
associated with $\lambda_1$. By theorem \ref{thm:bohm}
(proposition \ref{prop:a} and lemma \ref{lem:simple}), $\lambda_1$ is simple and, for $i$ large
enough,
\begin{equation}\label{lambda1}
\lambda_1<-2E,\qquad \lambda_2-\lambda_1>2E .
\end{equation}

The round $S^2$ factor also has $\Ric=Eg_{S^2}$, i.e. curvature $E$. Let $(x,y,z)$ be the coordinates of $\mathbb R^3\supset S^2$, and $Y$ the restriction to $S^2$
of the coordinate $z$, the zonal first spherical harmonic, so that $\Delta^{S^2}Y=2EY$. Let $\Delta_L^B$ be the
Lichnerowicz Laplacian of $g_B$. Since the connection and the curvature of the product split,
$\Delta_L^B(Yh_1)=(\Delta^{S^2}Y)h_1+Y\Delta_L^{B_1}h_1$; thus $Yh_1$ (with no $S^2$ component) is TT on $B$ and
$$
\Delta_L^B(Yh_1)=\lambda_*\,Yh_1,\qquad\lambda_*:=\lambda_1+2E<0 .
$$
$Yh_1$ is $H$-invariant and $\sigma$-odd. Let $\mathcal T$ be the space of $\hat G_1$-invariant TT tensors of $B_1$.

The next lemma identifies the bottom of the spectrum of the horizon on $H$-invariant tensors, and shows that the
eigenvalue $\lambda_*$ is carried by $Yh_1$ only.

\begin{lemma}\label{lem:bottom}
On $H$-invariant symmetric $2$-tensors of $B$, $\Delta_L^B\ge\lambda_1$, and the eigenspace for $\lambda_*$ is $\mathbb RYh_1$.
\end{lemma}
\begin{proof}
Write $h=h_{11}+h_{12}+h_{22}$ according to $T^*B=T^*B_1\oplus T^*S^2$: $h_{11}$ is a section of $S^2T^*B_1$, $h_{22}$ of
$S^2T^*S^2$, and $h_{12}$ of $T^*B_1\odot T^*S^2$ (all three depending on both factors). Since the connection and the
curvature of $g_B$ split, $\nabla^*\nabla$, $\Ric_B\circ$ (which is $E$ times the identity) and $\mathring R$ preserve the
three blocks, and $\mathring R$ only couples indices lying in the same factor; in particular $\mathring Rh_{12}=0$.
Hence $\Delta_L^B$ preserves the three blocks, and on mixed tensors
$\Delta_L^Bh_{12}=\nabla^*\nabla h_{12}+2Eh_{12}$, so its quadratic form is $\ge2E\|h_{12}\|^2$. On $h_{22}$ it acts as
$\Delta^{B_1}+\Delta_L^{S^2}$, where $\Delta^{B_1}$ is the rough Laplacian in the directions of $B_1$, acting on the
dependence on the first factor (the bundle $S^2T^*S^2$ being parallel along $B_1$, it acts on each component), and
$\Delta_L^{S^2}$ is the Lichnerowicz Laplacian of the second factor, acting on the dependence on it; similarly for
$h_{11}$ below. Now $\Delta_L^{S^2}\ge0$ on symmetric tensors of $S^2$ (they are of the form
$\varphi g+\mathring{\mathcal L}\omega$, and $\Delta_L$ commutes with $\mathring{\mathcal L}$ and maps $\varphi g$ to $(\Delta\varphi)g$).
On $h_{11}$ it acts as $\Delta_L^{B_1}+\Delta^{S^2}$. Decompose $h_{11}$, at each point of $S^2$, into its TT part for $B_1$ and the
rest (Lie derivatives and multiples of $g_{B_1}$, on which $\Delta_L^{B_1}\ge0$, because $\Delta_L\delta^*=\delta^*\Delta_H$
and $\Delta_L(\varphi g)=(\Delta\varphi)g$ on an Einstein manifold), and its dependence on the point of $S^2$ into spherical harmonics $Y_\ell$, with
$\Delta^{S^2}Y_\ell=\ell(\ell+1)EY_\ell$. These decompositions commute with $\Delta^B_L$. On the TT part in the harmonic
$\ell$ the eigenvalues are $\lambda_a+\ell(\ell+1)E$, where the $\lambda_a$ are the eigenvalues of $\Delta_L^{B_1}$ on
$\mathcal T$, all of them $\ge\lambda_1$ by the definition of $\lambda_1$; all the other contributions are $\ge0$. This proves $\Delta_L^B\ge\lambda_1$ (recall $\lambda_1<0$ by \eqref{lambda1}). Since
$\lambda_*<0$, the eigenvalue $\lambda_*$ can only come from $\lambda_a+\ell(\ell+1)E=\lambda_1+2E$ with
$\lambda_a\ge\lambda_1$, i.e. $\ell=1,\lambda_a=\lambda_1$, or $\ell=0,\lambda_a=\lambda_*$; the second case is excluded by
\eqref{lambda1}, since $\lambda_1<\lambda_*=\lambda_1+2E<\lambda_2$, so that $\lambda_*$ is not an eigenvalue on $\mathcal T$. The $K$-invariant first harmonics are the multiples of $Y$, and
$\lambda_1$ is simple.
\end{proof}

\section{The radial problem}\label{sec:radial}

We put $u:=1/r$; it is a defining function of the conformal infinity of $g_{r_h}$, and a tensor decays like $u^a$
if its $g_{r_h}$-norm is $O(u^a)$.

\emph{The decomposition into sectors.} The black hole \eqref{bh} is a warped product $g_{r_h}=g_2+r^2g_B$ of the
two-dimensional orbit space $(\mathbb R^2,g_2=f\,d\tau^2+f^{-1}dr^2)$ with $(B,g_B)$. Following Kodama and Ishibashi
\cite[\S2]{KI2} (see also \cite{KI}; their $n$ is our $N$), symmetric $2$-tensors on $M$ are expanded in harmonics
of $B$, which splits them into three types, or sectors:
\begin{itemize}
\item the \emph{tensor} type, $r^2\psi\,h_B$, with $\psi$ a function on the orbit space and $h_B$ a TT eigentensor of
$\Delta_L^B$;
\item the \emph{vector} type, built on a co-closed eigen-$1$-form $V$ of $B$: symmetric products of $V$ with
$1$-forms on the orbit space, and $\delta_B^*V$ times functions on the orbit space;
\item the \emph{scalar} type, built on an eigenfunction $S$ of $\Delta_B$: $S$ times symmetric
$2$-tensors on the orbit space, symmetric products of $dS$ with $1$-forms on the orbit space, and $Sg_B$ and
$(\nabla^2_BS)^\circ$ (the trace-free part of $\nabla^2_BS$) times functions on the orbit space.
\end{itemize}
Since $B$ is Einstein, the operators of $B$ which appear ($\delta_B$, $\delta_B^*$, $d$, $\Delta_L^B$, $\Delta_H^B$,
$\Delta_B$) intertwine these harmonics, so that $\Delta_L+2n$ and the linearised Einstein operator preserve each type
and each eigenspace on $B$, and the three types are $L^2$-orthogonal. An element of $\ker_{L^2}(\Delta_L+2n)$ is
therefore a sum of kernel elements of the three types, each built on a single eigentensor, eigenform or eigenfunction
of $B$, and an $H$-invariant kernel element involves only $H$-invariant harmonics. In the tensor sector the equation
reduces to a single radial ordinary differential equation (this section); in the vector and scalar sectors, the
linearised Einstein equations reduce, up to gauge, to master equations for one function on the orbit space
\cite{KI,KI2} (section \ref{sec:scalar}). Kernel elements are TT and decay like $u^n$ \cite{Lee}.

If $h_B$ is a TT tensor on $B$ with $\Delta_L^Bh_B=\lambda h_B$, a direct computation (the tensor modes of the warped
product \eqref{bh}, \cite{KI,GH}; see also \cite[\S4]{DelayTT}, where $\Delta_L$ is computed on the TT tensors of a
warped product with Einstein fibre over a one-dimensional base) gives, for $h=r^2\psi(r)h_B$ (TT for $g_{r_h}$),
\begin{equation}\label{radial}
(\Delta_L+2n)h=r^2\big(A_\lambda\psi\big)h_B,\qquad
A_\lambda\psi=-r^{-N}(r^Nf\psi')'+\frac{\lambda-2E}{r^2}\psi,
\end{equation}
$A_\lambda$ being self-adjoint on $L^2([r_h,\infty),r^Ndr)$ with quadratic form
$q_\lambda(\psi)=\int_{r_h}^\infty\big(r^Nf\psi'^2+(\lambda-2E)r^{N-2}\psi^2\big)dr$ (no boundary condition at $r_h$:
$r_h$ is the axis of the $\mathbb R^2$ factor, and the regular solutions are the smooth ones). Its essential spectrum is
$[\frac{(N+1)^2}4,\infty)$. Indeed, at the end $r=r_h$, $f$ has a simple zero and the solutions of $A_\lambda\psi=\mu\psi$
behave like $1$ and $\log(r-r_h)$: this end is limit circle and contributes nothing to the essential spectrum, the form
domain selecting the bounded solution. At infinity, with $r=e^\theta$ and $\psi=r^{-(N+1)/2}\tilde\psi$, which is unitary from
$L^2(r^Ndr)$ onto $L^2(d\theta)$, $A_\lambda$ becomes $-\frac{d^2\tilde\psi}{d\theta^2}+\big(\frac{(N+1)^2}4+O(e^{-2\theta})\big)\tilde\psi$ up
to terms whose coefficients are $O(e^{-2\theta})$, since $f=r^2(1+O(r^{-2}))$; the decomposition principle gives the
claim \cite{Weid}. Below the essential spectrum, the spectrum is discrete and consists of simple
eigenvalues (separated boundary conditions, \cite{Weid}). Let $\mu_0(r_h,\lambda)\le\mu_1(r_h,\lambda)\le\dots$ be the
min-max values of $q_\lambda$ with respect to the norm of $L^2(r^Ndr)$: those which lie below $\frac{(N+1)^2}4$ are these
eigenvalues, in increasing order, and the others are equal to $\frac{(N+1)^2}4$. Where they are eigenvalues, they are real-analytic in $r_h$ (after the change of variable $r=r_h\varrho$, which fixes the domain, the
operators form an analytic family of type (B) \cite[Ch.~VII]{Kato}), and strictly increasing in $\lambda$ (the potential
is, and the eigenfunctions do not vanish identically).

The following lemma shows that the lowest eigenvalue $\mu_0$ of $A_\lambda$ crosses zero, with a change of sign, when
$\lambda$ is below the Gibbons--Hartnoll threshold: it is negative for small black holes, by Hardy's inequality, and positive for large
ones, by comparison with the planar black brane.

\begin{lemma}\label{lem:cross}
Let $\lambda<\lambda_c:=\frac{(N-1)(9-N)}4$. Then $\mu_0(r_h,\lambda)<0$ for $r_h$ small and $\mu_0(r_h,\lambda)>0$ for
$r_h$ large. Hence $r_-(\lambda):=\sup\{r_h:\mu_0(r_h,\lambda)<0\}$ is finite and positive,
$\mu_0(r_-,\lambda)=0$, and there is $\iota>0$ such that $\mu_0(r_h,\lambda)<0$ for $r_h\in(r_--\iota,r_-)$ and
$\mu_0(r_h,\lambda)>0$ for $r_h\in(r_-,r_-+\iota)$.
\end{lemma}
\begin{proof}
\emph{Small $r_h$.} Since $E=N-1$, $\lambda_c-2E=-\frac{(N-1)^2}4$, so $2E-\lambda>\frac{(N-1)^2}4$. Let
$c\in(\frac{(N-1)^2}4,2E-\lambda)$. By lemma \ref{lem:hardy} (with $p=N$) there are $R>1$ and $\psi\in C^\infty_c((1,R))$ with
$\int r^N\psi'^2\le c\int r^{N-2}\psi^2$. For $\varsigma>0$ let $\psi_\varsigma(r)=\psi(r/\varsigma)$; both
integrals are multiplied by $\varsigma^{N-1}$. Since $f\le1+r^2$, for every $r_h<\varsigma$,
$$
\begin{aligned}
q_\lambda(\psi_\varsigma)&\le(1+\varsigma^2R^2)\int r^N\psi_\varsigma'^2+(\lambda-2E)\int r^{N-2}\psi_\varsigma^2\\
&\le \varsigma^{N-1}\big((1+\varsigma^2R^2)c+\lambda-2E\big)\int r^{N-2}\psi^2<0
\end{aligned}
$$
as soon as $(1+\varsigma^2R^2)c<2E-\lambda$. Hence $\mu_0(r_h,\lambda)<0$ for $r_h$ small.

\emph{Large $r_h$.} Put $r=r_h\varrho$. Then $f\ge r_h^2F(\varrho)$, $F(\varrho)=\varrho^2-\varrho^{1-N}$, and
$$
\frac{q_\lambda(\psi)}{\int r^N\psi^2}\ge\frac{\int_1^\infty\varrho^NF(\partial_\varrho\psi)^2\,d\varrho}{\int_1^\infty\varrho^N\psi^2d\varrho}
+\frac{\lambda-2E}{r_h^2}\ \ge\ \mu_{\rm pl}-\frac{|\lambda-2E|}{r_h^2},
$$
using $\varrho^{N-2}\le\varrho^N$, where $\mu_{\rm pl}$ is the bottom of the spectrum of the operator
$-\varrho^{-N}\partial_\varrho(\varrho^NF\partial_\varrho\psi)$ on $L^2([1,\infty),\varrho^Nd\varrho)$. This is the tensor-mode
operator, with zero potential, of the \emph{planar black brane}
$$
g_{\rm pl}=F(\varrho)d\tau^2+F(\varrho)^{-1}d\varrho^2+\varrho^2g_{\mathbb R^N},\qquad \varrho\ge1,
$$
which is a Poincar\'e--Einstein metric whose horizon $\{\varrho=1\}$ is the flat $\mathbb R^N$ (or a flat torus). Recall
that the metrics \eqref{bh} belong to a larger family \cite{Bir}: for $\mathsf K\in\{1,0,-1\}$ and an Einstein metric
$g_B$ with $\Ric_B=\mathsf K(N-1)g_B$, the metric $f_{\mathsf K}d\tau^2+f_{\mathsf K}^{-1}dr^2+r^2g_B$ with
$f_{\mathsf K}(r)=\mathsf K+r^2-\mathsf m\,r^{1-N}$, $\mathsf m$ a constant, is Einstein; \eqref{bh} is the case $\mathsf K=1$, and
$g_{\rm pl}$ is the case $\mathsf K=0$ with a flat horizon and $\mathsf m=1$. The planar black brane is also the limit of the
metrics $g_{r_h}$ as $r_h\to\infty$, written in the variables $\varrho=r/r_h$ and $r_h\tau$: indeed
$g_{r_h}=r_h^{-2}f\,d(r_h\tau)^2+(r_h^{-2}f)^{-1}d\varrho^2+\varrho^2\,r_h^2g_B$ with
$r_h^{-2}f=F(\varrho)+r_h^{-2}(1-\varrho^{1-N})$, and the metric $r_h^2g_B$ of the horizon, whose curvature tends to $0$, converges locally to the flat metric of
$\mathbb R^N$. The
quadratic form of this operator is $\ge0$; its essential spectrum is $[\frac{(N+1)^2}4,\infty)$ (same argument as for $A_\lambda$, since
$F$ has a simple zero at $\varrho=1$ and $F=\varrho^2(1+O(\varrho^{-N-1}))$); and $0$ is not an eigenvalue (the form vanishes
only on constants, which are not in $L^2$). Hence $\mu_{\rm pl}>0$, and $\mu_0>0$ for $r_h^2>|\lambda-2E|/\mu_{\rm pl}$.

\emph{Crossing.} By the above, the set $\{r_h:\mu_0(r_h,\lambda)<0\}$ is non-empty and bounded, so $r_-$ is finite and
positive. By definition of a supremum, $\mu_0(r_h,\lambda)\ge0$ for $r_h>r_-$, and there are $r_h<r_-$ arbitrarily close
to $r_-$ with $\mu_0(r_h,\lambda)<0$; by continuity, $\mu_0(r_-,\lambda)=0$. The function $\mu_0(\cdot,\lambda)$ is
real-analytic on the open set where it lies below the essential spectrum (simple isolated eigenvalue of an analytic
family, after the change of variable $r=r_h\varrho$ which fixes the domain), and this set contains its zeros. It is not
identically zero near $r_-$ (it takes negative values arbitrarily close to $r_-$, as just seen), so, being analytic,
it has an isolated zero at $r_-$: there
is $\iota>0$ such that $\mu_0(\cdot,\lambda)$ does not vanish on $(r_--\iota,r_-)\cup(r_-,r_-+\iota)$. It is
$\ge0$, hence $>0$, on $(r_-,r_-+\iota)$, and it is $<0$ on $(r_--\iota,r_-)$, an interval on which it does not
vanish and takes negative values.
\end{proof}

Let $r_-:=r_-(\lambda_*)$ (lemma \ref{lem:cross} applies since $\lambda_*<0<\lambda_c=3$ for $N=7$). For $r_h$ close
to $r_-$, the lowest min-max value $\mu_0(r_h,\lambda_*)$ is close to $\mu_0(r_-,\lambda_*)=0$, hence lies below the
essential spectrum $[\frac{(N+1)^2}4,\infty)$ and is therefore a simple eigenvalue of $A_{\lambda_*}$: we let
$\psi_{0,r_h}$ be a corresponding eigenfunction, normalised in $L^2([r_h,\infty),r^Ndr)$ and thus determined up to
sign, that is the ground state of $A_{\lambda_*}$ for $g_{r_h}$. We put $\psi_0:=\psi_{0,r_-}$ and
$$
k_{r_h}:=r^2\psi_{0,r_h}Yh_1,\qquad k_0:=k_{r_-},\qquad \nu(r_h):=\mu_0(r_h,\lambda_*).
$$
$k_{r_h}$ is an $H$-invariant, $\sigma$-odd, TT, $L^2$ eigentensor of $\Delta_L+2n$ for $g_{r_h}$ with eigenvalue
$\nu(r_h)$; $\nu$ is analytic near $r_-$, $\nu<0$ on $(r_--\iota,r_-)$, $\nu(r_-)=0$, $\nu>0$ on $(r_-,r_-+\iota)$.
Moreover $\nu(r_h)$ being a simple isolated eigenvalue of an analytic family, the associated spectral projection is
analytic in $r_h$, so that $r_h\mapsto\psi_{0,r_h}$, and hence $r_h\mapsto k_{r_h}$, are analytic in $L^2$ and,
by elliptic regularity, in the weighted H\"older spaces used in section \ref{sec:bif}.

The radius $r_-$ is the bifurcation radius. In order to exclude other radial zero modes at $r_-$, we need to know
that the second radial level of the lowest horizon eigenvalue is positive there. This follows from a comparison with
the planar black brane, since $r_-$ is large when $|\lambda_1|$ is large.

\begin{proposition}\label{prop:b}
There is $\Lambda_0>0$, depending only on $N$ and $E$, such that $\mu_1(r_-,\lambda_1)>0$ whenever $\lambda_1<-\Lambda_0$; moreover
$r_-^2\ge(|\lambda_1|-\mathsf b_*)/\mathsf b_*$ for a constant $\mathsf b_*>0$ depending only on $N$. In particular $r_-\ge6$ whenever
$\lambda_1<-\Lambda_1:=-37\mathsf b_*$.
\end{proposition}
\begin{proof}
The proof has five steps. Step 1 rescales the problem, so that the radial operator depends only on two parameters:
the size $\varepsilon=r_h^{-2}$ of the black hole and the coupling $\mathsf b$ of the attractive potential. Step 2
records three elementary properties of this two-parameter family. Step 3 studies the limiting planar operator
$\varepsilon=0$, which provides the two constants $\mathsf b_*$ and $\mu_*$ depending only on $N$. Step 4 locates
$r_-$ by comparing the two, and step 5 bounds the second eigenvalue at $r_-$.

\smallskip
\emph{Step 1: scaling.} Put $r=r_h\varrho$, $\varepsilon:=r_h^{-2}$. Then $f=r_h^2F_\varepsilon(\varrho)$ with
$$
F_\varepsilon(\varrho)=F(\varrho)+\varepsilon(1-\varrho^{1-N}),\qquad F(\varrho)=\varrho^2-\varrho^{1-N},
$$
and, for $\lambda<2E$, $A_\lambda$ becomes the operator
$$
P_{\varepsilon,\mathsf b}\psi=-\varrho^{-N}\partial_\varrho(\varrho^NF_\varepsilon\partial_\varrho\psi)-\frac{\mathsf b}{\varrho^2}\psi,\qquad
\mathsf b:=(2E-\lambda)\varepsilon>0 ,
$$
on $L^2([1,\infty),\varrho^Nd\varrho)$. Indeed, with $\partial_r=r_h^{-1}\partial_\varrho$ and $f=r_h^2F_\varepsilon$, the powers of $r_h$ cancel in
the kinetic term of $A_\lambda$, $r^{-N}\partial_r(r^Nf\partial_r\psi)=\varrho^{-N}\partial_\varrho(\varrho^NF_\varepsilon\partial_\varrho\psi)$, while the potential
$(\lambda-2E)r^{-2}$ becomes $-(2E-\lambda)r_h^{-2}\varrho^{-2}=-\mathsf b\varrho^{-2}$. Thus the operator only depends on the two
parameters $\varepsilon$ (the size of the black hole) and $\mathsf b$ (the strength of the attractive potential, which we
call the coupling). Let $\mu_j(\varepsilon,\mathsf b)$ be its min-max
values,
that is, its eigenvalues below the essential spectrum $[\frac{(N+1)^2}4,\infty)$, followed by
$\frac{(N+1)^2}4$; thus
$\mu_j(r_h,\lambda)=\mu_j(\varepsilon,\mathsf b)$.

\smallskip
\emph{Step 2: three elementary facts.} (i) Since $0\le1-\varrho^{1-N}\le\varrho^2-\varrho^{1-N}$ on $[1,\infty)$,
$F\le F_\varepsilon\le(1+\varepsilon)F$. The quadratic forms of the kinetic parts are therefore comparable, have the
same domain, and the min-max principle gives
\begin{equation}\label{sandwich}
\mu_j(0,\mathsf b)\le\mu_j(\varepsilon,\mathsf b)\le(1+\varepsilon)\,\mu_j\big(0,\tfrac{\mathsf b}{1+\varepsilon}\big).
\end{equation}
(ii) Since $0\le\varrho^{-2}\le1$, the multiplication by $\varrho^{-2}$ is bounded by $1$ on $L^2(\varrho^Nd\varrho)$, so each
$\mu_j(\varepsilon,\cdot)$ is non-increasing and $1$-Lipschitz; $\mu_0(\varepsilon,\cdot)$ is strictly decreasing where it is
an eigenvalue. (iii) The eigenvalues below the essential spectrum are simple: at $\varrho=1$, $F_\varepsilon$ has a simple
zero and the solutions behave like $1$ and $\log(\varrho-1)$, only the first one being in the form domain.

\smallskip
\emph{Step 3: the planar operator.} This is the case $\varepsilon=0$: $P_{0,\mathsf b}$ is the tensor-mode operator of the planar
black brane with the coupling $\mathsf b$ defined above. We have $\mu_0(0,0)=\mu_{\rm pl}>0$ (see the proof of lemma
\ref{lem:cross}), and $\mu_0(0,\mathsf b)\to-\infty$ as $\mathsf b\to\infty$: for a fixed $\psi\in C^\infty_c((1,\infty))$, the Rayleigh
quotient of $P_{0,\mathsf b}$ is an affine function of $\mathsf b$ with negative slope. By (ii), $\mu_0(0,\cdot)$ is continuous and
strictly decreasing where it is negative or zero, so there is a unique $\mathsf b_*>0$ with $\mu_0(0,\mathsf b_*)=0$. We define
$$
\mu_*:=\mu_1(0,\mathsf b_*)>0 .
$$
It is positive for the following reason: either $\mu_1(0,\mathsf b_*)$ is an eigenvalue below the essential spectrum, and
then it is strictly larger than the eigenvalue $\mu_0(0,\mathsf b_*)=0$ since the eigenvalues are simple by (iii); or it is
not, and then, by the definition of the min-max values, it equals the bottom $\frac{(N+1)^2}4>0$ of the essential
spectrum. Both $\mathsf b_*$ and $\mu_*$ only depend on $N$.

\smallskip
\emph{Step 4: location of $r_-$.} At $r_h=r_-$ and $\lambda=\lambda_*$ we have $2E-\lambda_*=|\lambda_1|$, hence
$\mathsf b_-:=|\lambda_1|\varepsilon_-$ with $\varepsilon_-:=r_-^{-2}$, and $\mu_0(\varepsilon_-,\mathsf b_-)=0$. By \eqref{sandwich},
$\mu_0(0,\mathsf b_-)\le0\le(1+\varepsilon_-)\mu_0(0,\mathsf b_-/(1+\varepsilon_-))$, so by the monotonicity (ii)
$$
\mathsf b_*\le \mathsf b_-\le(1+\varepsilon_-)\mathsf b_* .
$$
The second inequality reads $\varepsilon_-(|\lambda_1|-\mathsf b_*)\le \mathsf b_*$, so $\varepsilon_-\le \mathsf b_*/(|\lambda_1|-\mathsf b_*)$ as soon as
$|\lambda_1|>\mathsf b_*$.

\smallskip
\emph{Step 5: the second level.} We now consider the second radial eigenvalue $\mu_1$, at the same radius $r_-$, for the lowest
horizon eigenvalue $\lambda_1$ (the horizon mode $h_1$, without the factor $Y$); this is the most negative potential
among the $H$-invariant tensor modes, by lemma \ref{lem:bottom}, hence the smallest possible value of $\mu_1$. For
$\lambda=\lambda_1$ we have $2E-\lambda_1=|\lambda_1|+2E$, so the coupling at $r_-$ is
$$
\tilde{\mathsf b}_-:=(|\lambda_1|+2E)\varepsilon_-=\mathsf b_-\Big(1+\frac{2E}{|\lambda_1|}\Big)\ \ge\ \mathsf b_-\ \ge\ \mathsf b_* ,
$$
slightly larger than $\mathsf b_-$. We bound $\mu_1(r_-,\lambda_1)=\mu_1(\varepsilon_-,\tilde{\mathsf b}_-)$ from below in two steps. First,
by the left inequality of \eqref{sandwich}, $\mu_1(\varepsilon_-,\tilde{\mathsf b}_-)\ge\mu_1(0,\tilde{\mathsf b}_-)$. Second, since
$\mu_1(0,\cdot)$ is non-increasing and $1$-Lipschitz by (ii) and $\tilde{\mathsf b}_-\ge\mathsf b_*$,
$\mu_1(0,\tilde{\mathsf b}_-)\ge\mu_1(0,\mathsf b_*)-(\tilde{\mathsf b}_--\mathsf b_*)=\mu_*-(\tilde{\mathsf b}_--\mathsf b_*)$. Hence
$$
\mu_1(r_-,\lambda_1)\ge\mu_*-(\tilde{\mathsf b}_--\mathsf b_*) .
$$
It remains to see that $\tilde{\mathsf b}_--\mathsf b_*$ is small. By the location of $r_-$, $\mathsf b_-\le(1+\varepsilon_-)\mathsf b_*$, so that
$$
0\le \tilde{\mathsf b}_--\mathsf b_*\le \mathsf b_*\Big[(1+\varepsilon_-)\Big(1+\frac{2E}{|\lambda_1|}\Big)-1\Big],\qquad
\varepsilon_-\le\frac{\mathsf b_*}{|\lambda_1|-\mathsf b_*},
$$
and the right-hand side of the first inequality tends to $0$ as $|\lambda_1|\to\infty$. Hence there is a constant
$\Lambda_0$, depending only on $\mathsf b_*$, $\mu_*$ and $E$, i.e. on $N$ and $E$, such that $\tilde{\mathsf b}_--\mathsf b_*<\mu_*$, and
therefore $\mu_1(r_-,\lambda_1)>0$, whenever $|\lambda_1|>\Lambda_0$. Finally
$r_-^2=\varepsilon_-^{-1}\ge(|\lambda_1|-\mathsf b_*)/\mathsf b_*\ge36$ if $|\lambda_1|\ge37\mathsf b_*$.
\end{proof}

Since $\lambda_1^+(g_i)\to-\infty$, we may and do choose the B\"ohm metric $g_{B_1}$ so far in the sequence that, besides
\eqref{lambda1},
\begin{equation}\label{lambda1b}
\lambda_1<-\max(\Lambda_0,\Lambda_1),\qquad\text{so that}\qquad \mu_1(r_-,\lambda_1)>0\quad\text{and}\quad r_-\ge6 .
\end{equation}

\section{Vector and scalar perturbations}\label{sec:scalar}

\emph{Riemannian signature.} The $H$-invariant tensors are independent of $\tau$ and, by the invariance under
$\tau\mapsto-\tau$, have no components $d\tau\odot dr$ or $d\tau\odot\omega$, $\omega$ a $1$-form on $B$. For such tensors, the linearised Einstein
equations at $g_{r_h}$ are obtained from the static linearised equations at the Lorentzian metric
$-f\,d\tau_L^2+f^{-1}dr^2+r^2g_B$ by $\tau_L=i\tau$ (the static Einstein equations only involve $f=|\partial_{\tau_L}|^2$ up to its sign).
The Riemannian formulation of the master equations of \cite{KI}, the asymptotics of the master functions of $L^2$
kernel elements, and the reconstruction of a kernel element with vanishing master functions as a pure gauge tensor,
are carried out in \cite{CDK}; we use these results for the static modes, the base $B$ being Einstein instead of
maximally symmetric. For the master equations themselves this is not an issue: Kodama and Ishibashi observe that they
remain valid when the horizon is an arbitrary Einstein manifold \cite{KI}, and \cite{KI2} treats this case from the
outset (with a charge, which vanishes here). For the results of \cite{CDK}, it is legitimate because they use the base
only through a few identities satisfied by the harmonics, which hold on every compact Einstein manifold: these are
stated and proved in appendix \ref{app:harm}, where we also explain how they enter. Moreover $B$ has no scalar eigenvalue equal to the exceptional value $\kappa=N$, by
Lichnerowicz--Obata \cite{Lich58,Obata}, so that no exceptional scalar mode occurs.

We use the decomposition into sectors recalled in section \ref{sec:radial}; the tensor sector has been studied
there.

We first deal with the vector type and with the scalar modes built on constant functions, which can be handled
by symmetry, gauge and Birkhoff-type arguments.

\begin{lemma}\label{lem:vector}
\textup{(i)} There is no $H$-invariant kernel element of vector type. \textup{(ii)} The $H$-invariant kernel elements of
scalar type built on constant functions ($\kappa=0$) vanish unless $r_-^2=\frac{N-1}{N+1}$.
\end{lemma}
\begin{proof}
(i) We use the notation of \cite[\S5]{KI}: indices $a,b\in\{\tau,r\}$ refer to the orbit space with its metric
$f\,d\tau^2+f^{-1}dr^2$ and covariant derivative $D_a$, indices $i,j$ to $B$; $V$ is a co-closed $1$-form on $B$ with
$-D^jD_jV_i=k_V^2V_i$, and $V_{ij}=-\frac1{2k_V}(D_iV_j+D_jV_i)$ (if $k_V\ne0$). A vector-type tensor built on $V$ is
$h_{ai}=rf_aV_i$, $h_{ij}=2r^2H_TV_{ij}$, $h_{ab}=0$, with functions $f_a$, $H_T$ on the orbit space. Invariance under $\tau\mapsto-\tau$ gives $f_\tau=0$. If $V$ is not Killing, the gauge-invariant
quantity $F_a=f_a+\frac r{k_V}D_aH_T$ satisfies $D_a(r^{N-1}F^a)=0$ \cite[(4.6), (4.7b)]{KI2}, \cite[(5.10)--(5.12)]{KI}, hence, the perturbation being static,
$fF_r=C\,r^{1-N}$. The smoothness of $h$ along the axis $\{r=r_h\}$ (where $|dr|=\sqrt f\to0$) forces
$f_r=O(f^{-1/2})$, hence $C=0$ and $F=0$, i.e. $h=\mathcal L_Wg$ with $W$ proportional to $V$ and decaying. If $V$ is
Killing, $h=rf_r\,dr\odot V$ is also of the form $\mathcal L_Wg$, $W=\phi_V(r)V$ with $r^2\phi_V'=rf_r$, $\phi_V(\infty)=0$. In both
cases $h=2\delta^*W$ is TT, so $(\nabla^*\nabla+n)W=0$. The indicial roots of $\nabla^*\nabla+n$ on $1$-forms being
$-1$ and $n+1$ (for the pointwise norm of $g$, \cite{Lee}), the decaying solution $W$ is $O(u^{n+1})$ together with its
covariant derivative; hence the boundary terms in the integration by parts of $\langle(\nabla^*\nabla+n)W,W\rangle$
vanish, $\int|\nabla W|^2+n|W|^2=0$, and $W=0$.

(ii) Such a perturbation is constant along $B$, so that it satisfies the same equations as for a maximally
symmetric base (only $\Ric_B=(N-1)g_B$ enters). By the linearised Birkhoff theorem in all dimensions
\cite[\S3]{Kl} (see also \cite[App.~F]{CDK}, \cite[\S2]{KI}, \cite[\S5]{KI2}), it is, up to a gauge
transformation, the derivative of the family \eqref{bh} with respect to the mass, with the period $\tau_0$ adjusted for
regularity. It preserves the conformal infinity only if $d\tau_0/dr_h=0$, and
$\tau_0=4\pi r_h/((N+1)r_h^2+N-1)$ gives $r_h^2=\frac{N-1}{N+1}$; equivalently, the mass variation is in $L^2$ only
for this radius \cite[(2.20)--(2.21)]{CDK}.
\end{proof}

Scalar-type modes are built on $H$-invariant eigenfunctions $\Delta_BS=\kappa S$, $\kappa>0$.
By Lichnerowicz--Obata \cite{Lich58,Obata}, $\kappa>N$ ($B$ is not a round sphere); put $m:=\kappa-N>0$. By \cite[\S3, \S4]{KI} (the master equation also describes static perturbations) and \cite[\S3--5]{CDK}, the
master function $\Phi$ of an $H$-invariant kernel element of scalar type is a $\tau$-independent, hence radial, smooth
solution on the $(\tau,r)$-plane $\mathbb R^2$ of the master equation, which is $O(r^{-(N/2-1)})=O(r^{-5/2})$ at
infinity \cite[(3.7)]{CDK}; and if $\Phi=0$, the corresponding part of the kernel element is $\mathcal L_Wg$ with $W$
decaying \cite[\S4]{CDK}, hence vanishes, being TT. Thus such a kernel element is zero as soon as the static master equation
\begin{equation}\label{master}
-\frac{d^2\Phi}{dr_*^2}+V_S\Phi=0,\qquad \frac d{dr_*}=f\frac d{dr},\qquad V_S=\frac{f\,U}{16r^2\mathcal H^2},\quad \mathcal H=m+\tfrac{N(N+1)}2\mathsf x,
\end{equation}
$\mathsf x=r_h^{N-1}(1+r_h^2)r^{1-N}$, has no non-zero solution which is bounded at $r=r_h$ and $O(r^{-5/2})$ at infinity. Here
$U$ is the polynomial $Q(r)$ of \cite[(3.7)]{KI} for a horizon of positive curvature and a negative cosmological
constant (in their notation $K=1$ and $y=-r^2$), so that $f=1-\mathsf x+r^2$; for $N=7$,
$$
\begin{aligned}
U={}&\big(\mathfrak c_1\mathsf x^2-\mathfrak c_2m\mathsf x+\mathfrak c_3m^2\big)r^2+\mathfrak c_4\mathsf x^3+(\mathfrak c_5m+\mathfrak c_6)\mathsf x^2\\
&-(\mathfrak c_7m+\mathfrak c_8)m\mathsf x+16m^3+\mathfrak c_9m^2,
\end{aligned}
$$
with $\mathfrak c_1=197568$, $\mathfrak c_2=23520$, $\mathfrak c_3=60$, $\mathfrak c_4=153664$, $\mathfrak c_5=18144$, $\mathfrak c_6=47040$, $\mathfrak c_7=252$, $\mathfrak c_8=23520$, $\mathfrak c_9=252$.

For large eigenvalues $\kappa$, the potential of the master equation is positive, and a simple integration by
parts excludes zero modes.

\begin{proposition}\label{prop:scalarlarge}
For $N=7$, \eqref{master} has no such solution when $m\ge384\,(1+r_h^2)$.
\end{proposition}
\begin{proof}
On $[r_h,\infty)$, $0<\mathsf x\le \mathsf x(r_h)=1+r_h^2$, so $\mathsf r:=\mathsf x/m\le\frac1{384}$. The polynomial $\mathfrak c_1\mathsf r^2-\mathfrak c_2\mathsf r+\mathfrak c_3$ is decreasing on
$[0,\frac{\mathfrak c_2}{2\mathfrak c_1}]$, which contains $\frac1{384}$, and
$$
384^2\Big(\frac{\mathfrak c_1}{384^2}-\frac{\mathfrak c_2}{384}+\mathfrak c_3\Big)=197568-9031680+8847360=13248>0 ,
$$
so $\mathfrak c_1\mathsf x^2-\mathfrak c_2m\mathsf x+\mathfrak c_3m^2=m^2(\mathfrak c_1\mathsf r^2-\mathfrak c_2\mathsf r+\mathfrak c_3)>0$. Moreover
$$
16m^3-\mathfrak c_7m^2\mathsf x-\mathfrak c_8m\mathsf x=m^2(16m-\mathfrak c_7m\mathsf r-\mathfrak c_8\mathsf r)\ge m^2\Big(\big(16-\tfrac{21}{32}\big)m-\tfrac{245}4\Big)>0,
$$
since $\mathfrak c_7/384=\frac{21}{32}$, $\mathfrak c_8/384=\frac{245}4$ and $m\ge384$. All the other terms of $U$ are
non-negative, so $U>0$ and $V_S>0$. If $\Phi$ is a solution as above, multiplying \eqref{master} by $\Phi$ and integrating
in $r_*$ gives $\int(|d\Phi/dr_*|^2+V_S\Phi^2)dr_*=\big[\Phi\,d\Phi/dr_*\big]$. The boundary terms vanish. At the
horizon, $\Phi$ is a smooth radial function on the $(\tau,r)$-plane, that is, a smooth even function of the distance
$\mathsf d$ to the axis $\{r=r_h\}$; since $d\mathsf d/dr=f^{-1/2}$, we have
$d\Phi/dr_*=f\,d\Phi/dr=f^{1/2}\,d\Phi/d\mathsf d\to0$ as $r\to r_h$, while $\Phi$ has a limit. At infinity,
$\Phi=O(r^{-5/2})$ and $d\Phi/dr=O(r^{-7/2})$ (the solutions of \eqref{master} have expansions in powers of $r$ which can be
differentiated), so that $d\Phi/dr_*=f\,d\Phi/dr=O(r^{-3/2})$ and $\Phi\,d\Phi/dr_*=O(r^{-4})=O(u^4)$. Hence $\Phi=0$.
\end{proof}

\subsection*{An explicit positive supersolution}
We treat the remaining range $0<m<384(1+r_h^2)$ by the ground state transform (S-deformation). Recall that $u=1/r$, and put
$\mathsf s:=\mathsf xu^2$ and $\hat F:=fu^2=1+u^2-\mathsf s$, so that $d/dr_*=-\hat F\,d/du$, and note that $\hat F>0$ on
$(0,u_h)$, $u_h=1/r_h$. Let $\Psi$ be a smooth positive function of $r\in(r_h,\infty)$, or equivalently of
$u\in(0,u_h)$ (the supersolution), let $\omega:=d\log\Psi/dr_*$ be its logarithmic derivative, also a function of
$r$, and set $\tilde V:=V_S-d\omega/dr_*-\omega^2$. If $\Phi$ is a solution of \eqref{master} as in the statement above,
we define the function $\chi:=\Phi/\Psi$ on $(r_h,\infty)$; writing $\Phi=\Psi\chi$ and integrating by parts (appendix
\ref{app:pot}, A.1) gives
\begin{equation}\label{gst}
\int\Big(\Psi^2\Big|\frac{d\chi}{dr_*}\Big|^2+\tilde V\Phi^2\Big)dr_*=\Big[\Psi^2\chi\frac{d\chi}{dr_*}\Big]_{\rm horizon}^{\infty}.
\end{equation}

The following lemma reduces the absence of zero modes to the non-negativity of the modified potential $\tilde V$,
for an explicit choice of $\Psi$.

\begin{lemma}\label{lem:psi}
For $\varpi\ge0$ let $w:=\sqrt{9+(m+\varpi)u^2}$ and
$$
\Psi_\varpi:=\frac{u^{-1/2}\,e^{w}}{\mathcal H}\Big(\frac{w-3}{w+3}\Big)^{3/2},\qquad\text{so that}\qquad
\frac{d}{du}\log\Psi_\varpi=\frac{w-\frac12}u-\frac{\mathcal H_u}{\mathcal H} .
$$
Then $\Psi_\varpi>0$ on $(0,u_h]$, $\Psi_\varpi\sim c\,u^{5/2}$ as $u\to0$, and the boundary terms in \eqref{gst} vanish.
Consequently, if $\tilde V\ge0$ on $(0,u_h)$ and $\tilde V\not\equiv0$, \eqref{master} has no solution as above.
\end{lemma}
\begin{proof}
The formula for the derivative follows from $\frac d{du}\big(w+\frac32\log\frac{w-3}{w+3}\big)=\frac wu$. As $u\to0$,
$w-3\sim(m+\varpi)u^2/6$, hence $\Psi_\varpi\sim cu^{5/2}$, like $\Phi$, so that $\chi$ is bounded with bounded derivative
in $u$ and $\Psi^2\chi\,d\chi/dr_*=O(u^5)$. At the horizon $\Psi_\varpi$ and $\chi$ are smooth functions of $u$, while
$\frac{d\chi}{dr_*}=-\hat F\frac{d\chi}{du}\to0$, since $\hat F=fu^2$ vanishes there and $\frac{d\chi}{du}$ is
bounded; hence the boundary term $\Psi^2\chi\,\frac{d\chi}{dr_*}$ vanishes at the horizon as well.
If $\tilde V\ge0$, \eqref{gst} gives $d\chi/dr_*=0$ and $\tilde V\Phi^2=0$, so $\Phi=c\Psi_\varpi$ and $c=0$ where
$\tilde V>0$, hence everywhere.
\end{proof}

In appendix \ref{app:pot} we compute the modified potential $\tilde V$ of $\Psi_5$. It is convenient to use the
variables
$$
\vartheta:=w-3>0,\qquad \mathsf y:=m+\varpi=m+5,\qquad \mathsf s=1+u^2-\hat F\ge0 ,
$$
which are related by
\begin{equation}\label{tyu}
\vartheta(\vartheta+6)=w^2-9=\mathsf y\,u^2 .
\end{equation}
In these variables, the modified potential is given by an explicit polynomial with a simple structure.

\begin{lemma}\label{lem:G}
Let $N=7$ and $\varpi=5$. On $(0,u_h)$,
$$
\frac{16\mathcal H^2u^2w}{\hat F}\,\tilde V=\frac{\mathcal P(\vartheta,\mathsf y,\mathsf s)}{u^4\mathsf y^3},\qquad
\mathcal P=m^2\vartheta^3Q_0+\mathsf s\,\mathsf y\,m\,\vartheta^2Q_1+\mathsf s^2\mathsf y^2\vartheta\,Q_2+\mathsf s^3\mathsf y^3Q_3,
$$
where $m=\mathsf y-5$ and
\begin{align*}
Q_0&=16\,\mathsf y\,\vartheta(\vartheta+6)^2\\
&\quad+20736-65664\vartheta-70848\vartheta^2-26976\vartheta^3-4960\vartheta^4-448\vartheta^5-16\vartheta^6,\\
Q_1&=\mathsf y\,(13824+22656\vartheta+12032\vartheta^2+3552\vartheta^3+448\vartheta^4+16\vartheta^5)\\
&\quad+27648-296064\vartheta-276096\vartheta^2-101984\vartheta^3-16576\vartheta^4-976\vartheta^5,\\
Q_2&=64512\,m+\vartheta\,\tilde Q_2,\\
\tilde Q_2&=\mathsf y\,(98560+53760\vartheta+14336\vartheta^2+896\vartheta^3)\\
&\quad+96768-93184\vartheta-121856\vartheta^2-17024\vartheta^3,\\
Q_3&=12544\,\vartheta(\vartheta^2+4\vartheta+6).
\end{align*}
\end{lemma}
The proof is a direct computation, which we detail in appendix \ref{app:pot}; the factors $m^2$ and $m$ are not
imposed, they come out of the computation.

With this structure, the positivity of $\tilde V$ reduces to three one-variable inequalities.

\begin{proposition}[curved base]\label{prop:scalarK1}
Let $N=7$, $\varpi=5$ and $r_h\ge6$. Then for every $m$ with $0<m<384(1+r_h^2)$, $\tilde V>0$ on $(0,u_h)$. Hence
\eqref{master} has no solution as above.
\end{proposition}
\begin{proof}
On $(0,u_h)$ we have $\vartheta>0$, $\mathsf s\ge0$, $m>0$ and $\mathsf y>0$. Moreover:
\begin{enumerate}[label=(\roman*)]
\item since $u<u_h=1/r_h\le\frac16$, \eqref{tyu} gives $\mathsf y=\vartheta(\vartheta+6)/u^2>36\,\vartheta(\vartheta+6)$;
\item since $m<384(1+r_h^2)$, \eqref{tyu} gives $\vartheta(\vartheta+6)=\mathsf yu^2<(384(1+r_h^2)+5)/r_h^2<395$, hence $\vartheta<20$.
\end{enumerate}
The polynomials $Q_0$, $Q_1$ and $\tilde Q_2$ are affine in $\mathsf y$, with a coefficient of $\mathsf y$ which is positive for $\vartheta>0$. By (i)
they are therefore larger than their values at $\mathsf y=36\vartheta(\vartheta+6)$, which are (expanding $576\vartheta^2(\vartheta+6)^3$, etc.)
\begin{align*}
Q_0&>20736-65664\vartheta+53568\vartheta^2+35232\vartheta^3+5408\vartheta^4+128\vartheta^5-16\vartheta^6\\
&=1728\,(31\vartheta^2-38\vartheta+12)+16\vartheta^3(2202+338\vartheta+8\vartheta^2-\vartheta^3),\\
Q_1&>27648+2689920\vartheta+5115264\vartheta^2+3312544\vartheta^3\\
&\quad+1183808\vartheta^4+223664\vartheta^5+19584\vartheta^6+576\vartheta^7,\\
\tilde Q_2&>96768+21195776\vartheta+15038464\vartheta^2+5014912\vartheta^3\\
&\quad+709632\vartheta^4+32256\vartheta^5 .
\end{align*}
The discriminant of $31\vartheta^2-38\vartheta+12$ is $38^2-4\cdot12\cdot31=-44<0$, and by (ii), $8\vartheta^2-\vartheta^3>-12\vartheta^2>-240\vartheta$, so that
$2202+338\vartheta+8\vartheta^2-\vartheta^3>2202+98\vartheta>0$. Hence $Q_0,Q_1,\tilde Q_2>0$, so $Q_2=64512m+\vartheta \tilde Q_2>0$, and $Q_3>0$. Since $m,\vartheta,\mathsf y>0$ and $\mathsf s\ge0$,
$\mathcal P>0$, hence $\tilde V>0$ on $(0,u_h)$, and lemma \ref{lem:psi} concludes.
\end{proof}

Collecting the results of this section, we obtain:

\begin{corollary}\label{cor:c}
For $g_{r_-}$, $\ker_{L^2}(\Delta_L+2n)$ contains no $H$-invariant element of vector or scalar type.
\end{corollary}
\begin{proof}
Lemma \ref{lem:vector} treats the vector type and $\kappa=0$, since $r_-^2\ge36>\frac{N-1}{N+1}$ by \eqref{lambda1b}. For $\kappa>0$ we have
$m=\kappa-7>0$ (Lichnerowicz--Obata, $B$ not being a round sphere), and propositions \ref{prop:scalarlarge} and
\ref{prop:scalarK1} cover all $m>0$, since $r_-\ge6$.
\end{proof}

\section{The kernel at the bifurcation radius}\label{sec:kernel}

We can now put the previous sections together and determine the $H$-invariant kernel of $\Delta_L+2n$ at the
bifurcation radius.

\begin{proposition}\label{prop:kernel}
The $H$-invariant part of $\ker_{L^2}(\Delta_L+2n)$ for $g_{r_-}$ is $\mathbb Rk_0$. Consequently, for $r_h$ close to
$r_-$, $\nu(r_h)$ is the only eigenvalue of $\Delta_L+2n$ on $H$-invariant $L^2$ tensors close to $0$, and it is simple.
\end{proposition}
\begin{proof}
By corollary \ref{cor:c}, an $H$-invariant kernel element is a sum of tensor modes $r^2\psi(r)h_B$, $h_B$ an
$H$-invariant TT eigentensor of $\Delta_L^B$ with eigenvalue $\lambda$ and $A_\lambda\psi=0$ at $r_-$, i.e.
$\mu_j(r_-,\lambda)=0$ for some $j$. By lemma \ref{lem:bottom}, $\lambda\ge\lambda_1$, so for $j\ge1$,
$\mu_j(r_-,\lambda)\ge\mu_1(r_-,\lambda_1)>0$ by proposition \ref{prop:b} and the monotonicity in $\lambda$. Hence
$j=0$, and $\mu_0(r_-,\lambda)=0$ forces $\lambda=\lambda_*$, $\psi=c\psi_0$, and $h_B\in\mathbb RYh_1$ by lemma
\ref{lem:bottom}. The last statement follows from the continuity of the isolated eigenvalues of the analytic family
$\Delta_L+2n$ below its essential spectrum. The latter is contained in $[\frac{n^2}4,\infty)$: on trace-free
symmetric $2$-tensors of an asymptotically hyperbolic manifold of dimension $n+1$, the essential spectrum of
$\Delta_L$ is $[\frac{n(n-8)}4,\infty)$ \cite{DelayEss}, \cite[Prop.~D]{Lee}, i.e. that of $\Delta_L+2n$ is $[\frac{n^2}4,\infty)$; on
pure trace tensors $\Delta_L(\varphi g)=(\Delta\varphi)g$ and $\Delta+2n\ge2n$ has essential spectrum $[\frac{n^2}4+2n,\infty)$; and
restricting to the closed invariant subspace of $H$-invariant tensors does not enlarge the essential spectrum.
\end{proof}

\section{Bifurcation and non-extension}\label{sec:bif}

\begin{proof}[Proof of theorem \ref{thm:main}]
The proof has two parts. The first one constructs the bifurcating metrics in four steps: the linearised operator and
its kernel (step 1), the solution of the equation transverse to the kernel (step 2), the reduced equation and the
branch $v\mapsto g_v$ (step 3), and the regularity of these metrics (step 4). The second part shows that the rotation
about the $x$-axis does not extend to a Killing field of $g_v$, in four further steps (a)--(d), and then determines
which Killing fields of the $S^2$ factor do extend.

\smallskip
\emph{Step 1: the linearised operator and its kernel.} Fix $\frac n2<\zeta<n$, and let $\mathcal X$ be the space of $H$-invariant symmetric $2$-tensors in the weighted
H\"older space $C^{2,\alpha}_\zeta$ of \cite{Lee}. The metrics $g_{r_h}$ being uniformly equivalent in $C^2$ for $r_h$
near $r_-$, so are the corresponding weighted norms, and we regard $\mathcal X$ as a fixed Banach space. Put
$$
\mathcal E(r_h,k)=\Ric(g_{r_h}+k)+n(g_{r_h}+k)+\delta^*_{g_{r_h}+k}\beta_{g_{r_h}}(g_{r_h}+k)
$$
(Bianchi gauge relative to $g_{r_h}$, $\beta_g(k)=\delta_gk+\frac12d\,\mathrm{tr}_gk$). $\mathcal E$ is analytic, equivariant under
the isometries of $g_{r_h}$ (in particular under $\sigma$), $\mathcal E(r_h,0)=0$, and $D_k\mathcal E(r_h,0)=\frac12(\Delta_L+2n)$ is
Fredholm of index $0$ and $L^2$-symmetric \cite{Lee}. Its kernel on $\mathcal X$ at $r_h=r_-$ is
one-dimensional, spanned by $k_0$: since $\zeta>\frac n2$, the elements of this kernel are in $L^2$, and conversely the
elements of $\ker_{L^2}(\Delta_L+2n)$ decay like $u^n$ \cite{Lee}, hence lie in $C^{2,\alpha}_\zeta$ since $\zeta<n$; so this
kernel is the $H$-invariant part of $\ker_{L^2}(\Delta_L+2n)$, which is $\mathbb Rk_0$ by proposition \ref{prop:kernel}. By
symmetry and index $0$, the cokernel is also spanned by $k_0$. (Without $H$-invariance, the kernel would contain
$r^2\psi_0\tilde Yh_1$ for every first spherical harmonic $\tilde Y$ of $S^2$, and would have dimension at least $3$.)

\smallskip
\emph{Step 2: solving transverse to the kernel.} Let $P_{r_h}$ be the $L^2(g_{r_h})$-orthogonal projection onto
$\mathbb Rk_{r_h}$. By proposition \ref{prop:kernel}, $D_k\mathcal E(r_h,0)$ is an isomorphism from $(1-P_{r_h})\mathcal X$ onto its
image $(1-P_{r_h})\mathcal X^0$ (the corresponding space of $C^{0,\alpha}_\zeta$ tensors), uniformly for $r_h$ near $r_-$.
The implicit function theorem gives, for $v\in\mathbb R$ small, a unique small $k^\perp(r_h,v)\in(1-P_{r_h})\mathcal X$ with
$(1-P_{r_h})\mathcal E(r_h,vk_{r_h}+k^\perp)=0$, analytic, with $k^\perp(r_h,0)=0$. Moreover $\partial_vk^\perp(r_h,0)=0$:
differentiating the last equation in $v$ at $v=0$ gives
$(1-P_{r_h})D_k\mathcal E(r_h,0)\big(k_{r_h}+\partial_vk^\perp(r_h,0)\big)=0$; since
$D_k\mathcal E(r_h,0)k_{r_h}=\frac12\nu(r_h)k_{r_h}$ is in the range of $P_{r_h}$, this reads
$(1-P_{r_h})D_k\mathcal E(r_h,0)\partial_vk^\perp(r_h,0)=0$ with $\partial_vk^\perp(r_h,0)\in(1-P_{r_h})\mathcal X$, and the isomorphism
property gives $\partial_vk^\perp(r_h,0)=0$.

\smallskip
\emph{Step 3: the reduced equation and the branch.} It remains to solve the equation in the direction of the kernel.
Let
$$
\Xi(r_h,v):=\langle \mathcal E(r_h,vk_{r_h}+k^\perp(r_h,v)),k_{r_h}\rangle_{L^2(g_{r_h})}
$$
be the component of $\mathcal E$ along $k_{r_h}$, that is,
$P_{r_h}\mathcal E(r_h,vk_{r_h}+k^\perp(r_h,v))=\Xi(r_h,v)\,k_{r_h}/\|k_{r_h}\|^2$; $\Xi$ is the reduced (bifurcation) function,
a real-analytic function of two real variables. By construction, the small zeros of $\mathcal E(r_h,\cdot)$ in $\mathcal X$ are the
tensors $vk_{r_h}+k^\perp(r_h,v)$ with $\Xi(r_h,v)=0$. Its derivative with respect to $v$ at $v=0$ is, since
$\partial_vk^\perp(r_h,0)=0$,
$$
\partial_v\Xi(r_h,0)=\langle D_k\mathcal E(r_h,0)k_{r_h},k_{r_h}\rangle=\tfrac12\langle(\Delta_L+2n)k_{r_h},k_{r_h}\rangle
=\tfrac12\nu(r_h)\|k_{r_h}\|^2 .
$$
Since $\sigma k_{r_h}=-k_{r_h}$ and $\mathcal E$ is $\sigma$-equivariant, uniqueness
gives $k^\perp(r_h,-v)=\sigma k^\perp(r_h,v)$ and $\Xi(r_h,-v)=-\Xi(r_h,v)$. In particular $\Xi(r_h,0)=0$ (the trivial
branch $k=0$), so that $\Xi(r_h,v)=v\,\Theta(r_h,v)$ with $\Theta$ analytic,
even in $v$, and $\Theta(r_h,0)=\partial_v\Xi(r_h,0)=\frac12\nu(r_h)\|k_{r_h}\|^2$. The non-trivial solutions correspond to
the zeros of $\Theta$. Since $\nu$ changes sign at $r_-$
(lemma \ref{lem:cross}), for every small $\iota_1\in(0,\iota)$ and $v\ne0$ small enough,
$\Theta(r_--\iota_1,v)<0<\Theta(r_-+\iota_1,v)$, and there is $r(v)\in(r_--\iota_1,r_-+\iota_1)$ with
$\Theta(r(v),v)=0$ (intermediate value theorem); we take for $r(v)$ a zero of $\Theta(\cdot,v)$ in
$[r_--\iota,r_-+\iota]$ which is closest to $r_-$. Since $\Theta$ is even in $v$, we may choose $r(-v)=r(v)$, and since
$\iota_1$ is arbitrary, $r(v)\to r_-$ as $v\to0$. This gives zeros $k(v):=vk_{r(v)}+k^\perp(r(v),v)$ of $\mathcal E(r(v),\cdot)$,
and $\sigma k(v)=k(-v)$; moreover $k(v)-\sigma k(v)=2vk_{r(v)}+(k^\perp-\sigma k^\perp)(r(v),v)$, whose component along $k_{r(v)}$ is
$2vk_{r(v)}\ne0$, so $\sigma k(v)\ne k(v)$. As usual,
applying $\beta$ to $\mathcal E=0$ and using the injectivity of $\nabla^*\nabla+n$ on decaying $1$-forms, $\beta_{g_{r(v)}}(g_{r(v)}+k(v))=0$,
so $g_v:=g_{r(v)}+k(v)$ is Einstein; it is Poincar\'e--Einstein with the conformal infinity of $g_{r(v)}$ since $k(v)$
decays.

\smallskip
\emph{Step 4: regularity.} $k(v)$ is a $C^{2,\alpha}$ tensor on the whole of $M=\mathbb R^2\times B$ (including the
axis $\{r=r_{h}\}$, where $g_{r_h}$ is smooth), small in the weighted norm, so $g_v$ is a complete metric on $M$; it is
Einstein, hence smooth (indeed real-analytic) in the interior by elliptic regularity of the Einstein equation in
harmonic coordinates. At infinity, $g_v$ has the smooth conformal infinity $[d\tau^2\oplus g_B]$ ($\tau$ of period
$\tau_0(r(v))$), and by \cite{CDLS} it is polyhomogeneous. Since the boundary dimension $n=8$ is even, the
Fefferman--Graham expansion \cite{FG} shows that log terms can only appear through the obstruction tensor of the
conformal infinity, which is locally determined by it; it vanishes here because the same conformal infinity bounds the
smoothly compactifiable Poincar\'e--Einstein metric $g_{r(v)}$. Hence $g_v$ is smoothly conformally compact. Moreover
$k(v)=vk_0+o(v)$.

\smallskip
\emph{Non-extension of the rotation about the $x$-axis.} Let $X_0$ be the Killing field of $S^2$ generating the
rotations about the $x$-axis, extended to $M$ (it acts on the $S^2$ factor only); it is a Killing field of every
$g_{r_h}$, and its restriction to $\partial M$ is a Killing field of the conformal infinity. We argue by contradiction
and assume that, for a sequence of values of $v\ne0$ tending to $0$, $X_0|_{\partial M}$ extends to a Killing field $X_v$
of $g_v=g_{r(v)}+k(v)$, i.e. $X_v$ is a Killing field of $g_v$, smooth up to the boundary, which coincides with $X_0$
on $\partial M$. This argument has four steps.

\smallskip
\emph{(a) The tensor $\hat k$.} Put $\hat k:=\mathcal L_{X_0}k_0=r^2\psi_0\,(X_0Y)\,h_1$. Since $X_0$ is a Killing field of
$g_{r_-}$, $\mathcal L_{X_0}$ commutes with $\Delta_L$, with the trace and with the divergence of $g_{r_-}$, so that $\hat k$ is
a TT element of $\ker_{L^2}(\Delta_L+2n)$ for $g_{r_-}$, decaying like $u^n$. It is not zero, since $X_0Y=\pm y\not\equiv0$
($Y=z$ and $X_0=\pm(y\partial_z-z\partial_y)$). Note that $\hat k$ is not $H$-invariant, which is consistent with proposition
\ref{prop:kernel}. Moreover $\mathcal L_{X_0}\hat k=r^2\psi_0(X_0^2Y)h_1=-k_0$, and, the flow of $X_0$ being isometric,
$\mathcal L_{X_0}$ is skew-symmetric for the $L^2(g_{r_-})$ product: for every tensor $k$ in $C^{2,\alpha}_\zeta$,
\begin{equation}\label{skew}
\langle\mathcal L_{X_0}k,\hat k\rangle_{L^2(g_{r_-})}=-\langle k,\mathcal L_{X_0}\hat k\rangle_{L^2(g_{r_-})}=\langle k,k_0\rangle_{L^2(g_{r_-})} .
\end{equation}

\smallskip
\emph{(b) The difference $Z_v:=X_v-X_0$ is $O(v)$.} A Killing field $X$ of an Einstein metric $g$ with $\Ric=-ng$ satisfies
$(\nabla^*\nabla+n)X^\flat=0$ (apply the Bianchi operator $\beta$ to $\delta^*X^\flat=0$, and use
$2\beta\delta^*=\nabla^*\nabla-\Ric$). Thus $X_v$ satisfies this equation for $g_v$, and $X_0$ for $g_{r(v)}$. Hence
$$
(\nabla^*\nabla+n)_{g_v}Z_v^{\flat_v}=-(\nabla^*\nabla+n)_{g_v}X_0^{\flat_v}+(\nabla^*\nabla+n)_{g_{r(v)}}X_0^{\flat_{r(v)}}=:\mathsf R_v ,
$$
where $\flat_v$ and $\flat_{r(v)}$ denote the musical isomorphisms of $g_v$ and of $g_{r(v)}$. Two things change
between the two terms. First the operator: their difference is a second order operator whose coefficients are
bounded by $k(v)$ and its first two covariant derivatives, which are $O(v\,u^\zeta)$, while $X_0$ and its covariant
derivatives are $O(u^{-1})$ (it is the lift of a Killing field of the boundary); this contributes $O(v)$ in
$C^{0,\alpha}_{\zeta-1}$. Then the $1$-form, with
$X_0^{\flat_v}-X_0^{\flat_{r(v)}}=k(v)(X_0,\cdot)=O(v\,u^{\zeta-1})$, hence $O(v)$ in $C^{2,\alpha}_{\zeta-1}$ together
with its first two derivatives; since $\nabla^*\nabla+n$ preserves the weight, this contributes $O(v)$ in
$C^{0,\alpha}_{\zeta-1}$ as well. Therefore $\mathsf R_v=O(v)$ in $C^{0,\alpha}_{\zeta-1}$. On the other hand $X_v$ and $X_0$ are smooth up to
the boundary and coincide there, so that $Z_v$ vanishes on $\partial M$ and $|Z_v|_{g_v}=O(1)$. Since
$\nabla^*\nabla+n$ is an isomorphism on $1$-forms between the weighted spaces $C^{2,\alpha}_s$ and
$C^{0,\alpha}_s$ for $-1<s<n+1$ \cite{Lee} (its indicial roots are $-1$ and $n+1$, and it has no $L^2$ kernel),
uniformly for the metrics $g_v$, $Z_v^{\flat_v}$ is the unique solution in the space of weight $s=0$, and it coincides
with the solution in the space of weight $\zeta-1\in(0,n+1)$ (the norms below being those of $g_v$, for which
$Z_v$ and $Z_v^{\flat_v}$ have the same size):
$$
\|Z_v\|_{C^{2,\alpha}_{\zeta-1}}\le C|v| .
$$

\smallskip
\emph{(c) The identity.} From $\mathcal L_{X_v}g_v=0$, $\mathcal L_{X_0}g_{r(v)}=0$ and $X_v=X_0+Z_v$,
$$
0=\mathcal L_{X_0}\big(g_{r(v)}+k(v)\big)+\mathcal L_{Z_v}g_v=\mathcal L_{X_0}k(v)+\mathcal L_{Z_v}g_v,
$$
that is,
$$
\mathcal L_{X_0}k(v)=-\mathcal L_{Z_v}g_v .
$$

\smallskip
\emph{(d) Pairing with $\hat k$.} We take the $L^2(g_{r_-})$ product of both sides with $\hat k$ and compare the terms of
order $v$. On the left, by \eqref{skew} and $k(v)=vk_0+o(v)$ in $C^{2,\alpha}_\zeta$,
$$
\langle\mathcal L_{X_0}k(v),\hat k\rangle=\langle k(v),k_0\rangle=v\|k_0\|^2+o(v) .
$$
On the right, write $\mathcal L_{Z_v}g_v=\mathcal L_{Z_v}g_{r_-}+\mathcal L_{Z_v}(g_v-g_{r_-})$. Since $g_v\to g_{r_-}$
in $C^2$ uniformly on $M$, the weighted norms of $g_v$ and $g_{r_-}$ are uniformly equivalent, so the bound of
\emph{(b)} holds for those of $g_{r_-}$ as well, which we use from now on. The first term is
$2\delta^*_{g_{r_-}}Z_v^{\flat_{r_-}}$ (the dual being that of $g_{r_-}$, the metric which is differentiated), and an
integration by parts on $\{u\ge u_0\}$ gives
$$
\langle\mathcal L_{Z_v}g_{r_-},\hat k\rangle=2\langle Z_v^{\flat_{r_-}},\div\hat k\rangle+\lim_{u_0\to0}(\text{boundary term})=0,
$$
because $\hat k$ is divergence-free and the boundary term on $\{u=u_0\}$ is bounded by
$\int_{\{u=u_0\}}|Z_v||\hat k|=O(u_0^{\zeta-1}u_0^{n}u_0^{-n})=O(u_0^{\zeta-1})\to0$. For the second term,
$g_v-g_{r_-}=(g_{r(v)}-g_{r_-})+k(v)$ tends to $0$ with its covariant derivative, uniformly on $M$, as $v\to0$ (recall
that $r(v)\to r_-$), so that
$|\mathcal L_{Z_v}(g_v-g_{r_-})|\le C\big(|Z_v||\nabla(g_v-g_{r_-})|+|\nabla Z_v||g_v-g_{r_-}|\big)=o(1)\cdot O(v\,u^{\zeta-1})$;
since $|\hat k|=O(u^n)$ and the volume form is $O(u^{-n-1})du$, the product is integrable ($\zeta>1$) and
$\langle\mathcal L_{Z_v}(g_v-g_{r_-}),\hat k\rangle=o(v)$. Altogether
$$
v\|k_0\|^2+o(v)=o(v),
$$
which is impossible for $v\ne0$ small. Hence, for $0<|v|$ small, the rotation $X_0$ about the $x$-axis does not extend
to a Killing field of $g_v$.

\smallskip
\emph{All the Killing fields of $S^2$.} We finally determine which Killing fields of the $S^2$ factor extend. Let
$X_x,X_y,X_z$ be the generators of the rotations about the three axes; the Killing fields of $S^2$ are the
$X=c_xX_x+c_yX_y+c_zX_z$, $c_x,c_y,c_z\in\mathbb R$, and the Lie algebra of $K$ is $\mathbb RX_z$. We claim that
\emph{$X|_{\partial M}$ extends to a Killing field of $g_v$ if and only if $c_x=c_y=0$}, which is the last statement of
theorem \ref{thm:main}. Indeed, $g_v$ is $K$-invariant, so $X_z$ is a Killing field of $g_v$; hence $X$ extends if and
only if $X^\perp:=c_xX_x+c_yX_y$ does. If $(c_x,c_y)\ne(0,0)$, $X^\perp$ is the image of $(c_x^2+c_y^2)^{1/2}X_x$ under a
rotation about the $z$-axis, that is, under an element of $K$. Since the elements of $K$ are isometries of $g_v$, the
field $X^\perp|_{\partial M}$ extends to a Killing field of $g_v$ if and only if $X_x|_{\partial M}$ does; but $X_x=X_0$
does not extend, by the above. Hence $X^\perp$, and therefore $X$, does not extend.

Finally, the hypotheses of the extension problem are satisfied: \eqref{pi1} holds, and the conformal infinity
$[d\tau^2\oplus g_{B_1}\oplus g_{S^2}]$ is the product of a circle with the Einstein metric $g_B$, for which the rotations of
$S^2$ are Killing fields.
\end{proof}

\begin{proof}[Proof of theorems \ref{thm:bohm} and \ref{thm:scalar}]
Theorem \ref{thm:bohm} is the content of propositions \ref{prop:model}, \ref{prop:a} and of lemma
\ref{lem:simple}. Theorem \ref{thm:scalar} is propositions \ref{prop:scalarlarge} and \ref{prop:scalarK1}, together with the reduction explained
in section \ref{sec:scalar} (for $\kappa>7$, $m=\kappa-7>0$; by Lichnerowicz--Obata \cite{Lich58,Obata}, $\kappa>7$
for every $\kappa>0$ when $B$ is not a round sphere).
\end{proof}

\section{Remarks on the extension of Killing fields}\label{sec:anderson}

We conclude with four remarks relating theorem \ref{thm:main} to \cite{An2}, to dimension four and to the Lorentzian
setting. The first one indicates the step of the proof of \cite[Th.~5.4]{An2} which, in our reading, does not apply in
general.

\begin{remark}[On the proof of {\cite[Th.~5.4]{An2}}]\label{rem:gap}
We use the notation of the proof of \cite[Th.~5.4]{An2}, in \cite[\S5]{An2} (in this remark and the next one, $t$, $\tau$, $A$, $\kappa$, and $h$, $h_1$, $h_2$ in the quotation, are not the objects of the previous sections): $t$ is a geodesic defining function, $\mathbf N=-t\partial_t$,
$\tau=A-\mathsf Hg_t$ on the level set $S(t)$ ($A$ its second fundamental form, $\mathsf H$ its mean curvature), $X$ is the
extension of the boundary Killing field in Bianchi gauge, $\kappa=\delta^*X$ and $\tilde\kappa=t^{-n}\kappa$. In \cite{An2} a prime denotes
the linearisation, $\tau'_h=D\tau(h)$, and $\hat E$ is the gauged Einstein operator; we quote the two formulas
verbatim (the reference \cite[(5.23)--(5.24)]{An2} is to the corresponding identities in the compact case):

\smallskip
\begin{quote}
``Now as in (5.23)-(5.24) (setting $M=B(t)$) one then has
\begin{equation*}
\tag{5.25}
\frac12\int_{S(t)}\langle\mathcal L_X\tau,h\rangle=\int_{S(t)}\langle\tau'_h,(\delta^*X)^T\rangle
+\int_{B(t)}\langle\hat E'_h,\delta^*X\rangle=\int_{S(t)}\langle\tau'_h,(\delta^*X)^T\rangle .
\end{equation*}
The equation (5.25) holds for any $h$ such that $\tilde L(h)=0$ satisfying the boundary conditions $\delta(h)=0$,
$[h^T]_0=h_1$, $H'_h=h_2$ with $h_1$, $h_2$ arbitrary on $S(t)$. [\dots] As in [2] [that is, \cite{AH}], let $\kappa=\delta^*X$ and
$\tilde\kappa=t^{-n}\delta^*X$, so that $\tilde\kappa$ is uniformly bounded as $t\to0$; [\dots]
We then choose $h$ such that $[h^T]_0=[\tilde\kappa^T]_0$ on $S(t)$. Substituting this in (5.25) gives the basic
relation
\begin{equation*}
\tag{5.26}
\int_{S(t)}\langle\mathcal L_X\tau,\tilde\kappa\rangle=2\int_{S(t)}\langle\tau'_{\tilde\kappa},(\delta^*X)^T\rangle ."
\end{equation*}
\end{quote}
\smallskip

\noindent The proof then compares two computations of the two sides of (5.26).
The left side is computed as $-\frac{3n-2}4\int_{\partial M}|\mathcal L_Xg_{(n)}|^2+o(1)$ and the right side as
$-\frac{2n-2}4\int_{\partial M}|\mathcal L_Xg_{(n)}|^2+o(1)$, whence $\mathcal L_Xg_{(n)}=0$.

In passing from (5.25) to (5.26), $D\tau(h)$ is replaced by $D\tau(\tilde\kappa)$. But (5.25) is established
for a solution $h$ of $\tilde L(h)=0$ on $B(t)$, and $D\tau(h)$ depends on the normal derivative of $h$ at $S(t)$, which
is not determined by the boundary condition $[h^T]_0=[\tilde\kappa^T]_0$. On $S(t)$, $\tilde\kappa$ coincides with the
constant multiple $t^{-n}\kappa$, so that the natural solution with these boundary values is $h=c\,\kappa$, $c$
constant. Since
$\mathcal L_Xg=2\kappa$ and $X$ is tangent to $S(t)$ up to negligible terms, $\mathcal L_X\tau=D\tau(2\kappa)=2D\tau(\kappa)$, and
\textup{(5.26)} then reduces to an identity. More precisely, with $\kappa^T=\frac12t^{n-2}\mathcal L_Xg_{(n)}+\dots$ one
has
$$
\begin{gathered}
\nabla_{\mathbf N}\kappa=-n\kappa+\dots,\qquad DA(\kappa)=\tfrac12\nabla_{\mathbf N}\kappa+\kappa=-\tfrac{n-2}2\kappa,\\
D\tau(\kappa)=DA(\kappa)-n\kappa=-\tfrac{3n-2}2\kappa
\end{gathered}
$$
(the trace terms are negligible). In our reading, both sides are therefore equal to
$-\frac{3n-2}4\int|\mathcal L_Xg_{(n)}|^2$, and no information on $\mathcal L_Xg_{(n)}$ follows. The value
$-\frac{2n-2}4$ is obtained by applying the variation formula to the field $t^{-n}\kappa$ itself, for which
$\nabla_{\mathbf N}(t^{-n}\kappa)=t^{-n}(\nabla_{\mathbf N}\kappa+n\kappa)=O(t)$ instead of $-n\,t^{-n}\kappa$; the difference
$n\cdot\frac14\int|\mathcal L_Xg_{(n)}|^2$ is exactly the discrepancy between the two values. Since $t^{-n}\kappa$ is not
a solution of the linearised equations, \textup{(5.25)} does not seem to apply to it.

\end{remark}

The second one discusses which weakenings of the hypothesis $\mathcal L_{X_0}g_{(n)}=0$, under which the extension is
known \cite[Cor.~4.4]{AH} (see also \cite[\S5.2]{CD11}; in Lorentzian signature, see \cite[Th.~1.6]{CD11} for stationary metrics, and \cite{HS} for general vacuum metrics, near boundary domains satisfying a convexity condition), are possible.

\begin{remark}[Weaker hypotheses]\label{rem:weaker}
Let $(M,g)$ be a smoothly conformally compact Poincar\'e--Einstein metric with $\partial M$ connected, $X_0$ a Killing
field of a representative $\gamma$ of its conformal infinity, $t$ the geodesic defining function associated with
$\gamma$, $g=t^{-2}(dt^2+g_t)$ and $g_t=\gamma+\dots+t^ng_{(n)}+\dots$. The flow of $X_0$ preserves $\gamma$, hence $t$
and the locally determined coefficients of the expansion, so that $\mathcal L_{X_0}g_{(n)}$ is the only quantity
which is not automatically zero.

\emph{(a) A trace condition gives nothing.} The trace $\mathrm{tr}_\gamma g_{(n)}$ and the divergence
$\div_\gamma g_{(n)}$ are given by local expressions in $\gamma$ \cite{FG} (both vanish when $n$ is odd), hence are
invariant under $X_0$, and since $\mathcal L_{X_0}\gamma=0$,
$$
\mathrm{tr}_\gamma(\mathcal L_{X_0}g_{(n)})=X_0(\mathrm{tr}_\gamma g_{(n)})=0,\qquad
\div_\gamma(\mathcal L_{X_0}g_{(n)})=\mathcal L_{X_0}(\div_\gamma g_{(n)})=0 .
$$
Thus $\mathcal L_{X_0}g_{(n)}$ is always a TT tensor on $(\partial M,\gamma)$; in theorem \ref{thm:main} it is a
non-zero one. These identities are the limits at infinity of the constraint equations. For a boundary at a finite
distance, with second fundamental form $A$ and mean curvature $\mathsf H$, the Codazzi and Gauss equations
$\div_\gamma A=-d\mathsf H$ and $R_\gamma-|A|^2_\gamma+\mathsf H^2=-n(n-1)$ give, for a Killing field $X_0$ of $\gamma$,
$$
\begin{gathered}
\mathrm{tr}_\gamma(\mathcal L_{X_0}A)=X_0(\mathsf H),\qquad \div_\gamma(\mathcal L_{X_0}A)=-d\big(X_0(\mathsf H)\big),\\
\langle A,\mathcal L_{X_0}A\rangle_\gamma=\mathsf H\,X_0(\mathsf H).
\end{gathered}
$$
There the condition $X_0(\mathsf H)=0$, from which the strategy of \cite{An2} starts, is a genuine hypothesis (and a Killing field of $\gamma$ which also preserves $A$
extends to a Killing field near the boundary \cite[Th.~5.1]{CD11}); at infinity its analogue is automatic.

\emph{(b) Localisation.} If $\mathcal L_{X_0}g_{(n)}$ vanishes on a non-empty open set $\mathcal U\subset\partial M$, it
vanishes on $\partial M$. Indeed, let $\mathcal O$ be the interior of the zero set of $\mathcal L_{X_0}g_{(n)}$, let $X$ be
the solution of the Killing extension problem $(\nabla^*\nabla+n)X^\flat=0$ with boundary value $X_0$
\textup{(}section \ref{sec:intro}\textup{)} and $\kappa=\delta^*X^\flat$, an infinitesimal Einstein deformation with
$\kappa^T=\frac12t^{n-2}\mathcal L_{X_0}g_{(n)}+\dots$ (remark \ref{rem:gap}). Let $o\in\overline{\mathcal O}$ and let
$\mathcal W_1\subset\partial M$ be a ball around $o$; $\mathcal W_0:=\mathcal W_1\cap\mathcal O$ is open and non-empty. The Cauchy data of $\kappa$
vanish on $\mathcal W_0$, so that, by unique continuation from the conformal infinity \cite{Bi,AH}, \cite[\S5.2]{CD11}, $\kappa=\delta^*W^\flat$ near $\mathcal W_0$
with $W=0$ on $\mathcal W_0$, and $Z:=X-W$ is a Killing field near $\mathcal W_0$ with boundary value $X_0$. Since $g$ is real-analytic
in the interior, $Z$ extends to a Killing field on the simply connected set $\mathcal W_1\times(0,t_0)$ \cite{Nomizu}. Its
boundary value is a conformal Killing field on $\mathcal W_1$ which coincides with $X_0$ on $\mathcal W_0$, hence on $\mathcal W_1$; its flow
preserves $g$ and acts on the boundary by isometries of $\gamma$, hence preserves $t$ and the expansion of $g_t$, and
$\mathcal L_{X_0}g_{(n)}=0$ on $\mathcal W_1$. Thus $o\in\mathcal O$, and $\mathcal O$, open, closed and non-empty, is $\partial M$.

If moreover $\pi_1(M,\partial M)=0$, then $X_0$ extends to a Killing field of $g$. Indeed, the Cauchy data of
$\kappa$ now vanish on $\partial M$, so $\kappa=\delta^*W^\flat$ near $\partial M$ with $W=0$ on $\partial M$ \cite{AH},
and $Z=X-W$ is a Killing field on a collar $\mathcal C$ of $\partial M$. Since $\pi_1(\mathcal C)=\pi_1(\partial M)\to\pi_1(M)$ is onto,
every loop of $M$ based in $\mathcal C$ is homotopic to a loop in $\mathcal C$, along which the analytic continuation of $Z$ is
trivial; hence $Z$ extends to a Killing field of $(M,g)$ \cite{Nomizu}, with boundary value $X_0$. This is the route of
\cite[Cor.~4.4]{AH}; without the topological hypothesis, the extension near $\partial M$ given by \cite{AH}, \cite[\S5.2]{CD11} need not be global \cite[Rem.~5.2]{CD11}. In particular, for the metrics $g_v$ of theorem \ref{thm:main}, which satisfy \eqref{pi1},
$\mathcal L_{X_0}g_{(n)}$ vanishes on no non-empty open subset of $\partial M$.
\end{remark}

The construction also applies in Lorentzian signature.

\begin{remark}[Lorentzian signature]\label{rem:lorentz}
The metrics $g_v$ are static: being invariant under the translations and the reflection $\tau\mapsto-\tau$, they read
$g_v=\mathsf V^2d\tau^2+\hat g$ away from the axis $\{\mathsf V=0\}$, with $\mathsf V$ and $\hat g$ independent of $\tau$. The Einstein
equations $\Ric=-ng$ for such a metric read $\Ric(\hat g)-\mathsf V^{-1}\hat\nabla^2\mathsf V=-n\hat g$, $\mathrm{tr}_{\hat g}\hat\nabla^2\mathsf V=n\mathsf V$
in both signatures, so that $g^L_v:=-\mathsf V^2d\tau_L^2+\hat g$ is a static vacuum metric with negative cosmological constant,
asymptotically anti-de Sitter with conformal boundary $(\mathbb R\times B,[-d\tau_L^2\oplus g_B])$ (the exterior of a black
hole whose bifurcation surface is the axis). The Killing fields $X$ of $S^2$ which are not in the Lie algebra of $K$
are Killing for $-d\tau_L^2\oplus g_B$, but do not extend to Killing fields of $g^L_v$, even near the conformal boundary.
Indeed, a Killing field of $g^L_v$ whose boundary value vanishes is zero near the boundary (its flow fixes the
boundary pointwise and preserves the geodesic defining function). If a Killing field $\tilde X$ of $g^L_v$ extended $X$, then $[\partial_{\tau_L},\tilde X]$ and
$\tilde X$ minus its image under $\tau_L\mapsto-\tau_L$ would be such fields, so that $\tilde X$ would be independent of $\tau_L$ and tangent to
the slices; $\mathcal L_{\tilde X}g^L_v=0$ then means $\mathcal L_{\tilde X}\hat g=0$ and $\tilde X(\mathsf V)=0$, i.e. $\tilde X$ would be a Killing field of
$g_v$ near $\partial M$ with boundary value $X$, which forces $\mathcal L_Xg_{(n)}=0$ there (see remark
\ref{rem:weaker}(b)), a contradiction. This is consistent with \cite[Th.~1.6]{CD11} and \cite{HS}, since $X$ does not
preserve the undetermined term; and the exterior region, diffeomorphic to $\mathbb R\times(r_h,\infty)\times B$, is
simply connected.
\end{remark}

The last one concerns dimension four, where \cite[Th.~5.4]{An2} has been used.

\begin{remark}[Dimension four]\label{rem:dim4}
The step discussed in remark \ref{rem:gap} uses neither the dimension nor the topology, so \cite[Th.~5.4]{An2} is
equally open in dimension four. The statement is used in \cite[Prop.~1]{ACD}, for a $C^{m,\alpha}$ conformal
compactification with $m\ge4$, to show that every asymptotically hyperbolic Einstein metric filling
in a strictly globally static boundary metric carries an isometric $S^1$ action, hence is static. \textup{(}The version
quoted in \cite{ACD} carries no topological hypothesis, but the manifolds considered there, $S^1\times\mathbb R^3$,
$\mathbb R^2\times S^2$ and $D^2\times T^2$, all satisfy $\pi_1(M,\partial M)=0$.\textup{)}

In \cite{An3}, this extension statement is \cite[Th.~2.5]{An3}, which is attributed to \cite{AH}; by remark \ref{rem:weaker}, the extension is known under the
additional hypothesis $\mathcal L_{X_0}g_{(n)}=0$. It is used in \cite{An3} for the computations of the degree of the boundary map in
\cite[\S7]{An3}, in particular in the proof of \cite[Th.~B]{An3} and in \cite[Prop.~7.2]{An3}, where it is also used to
show that the $L^2$ kernel along the anti-de Sitter--Schwarzschild metrics inherits their symmetries.

We know of no counterexample in dimension four, and the mechanism of this paper is not available there. The horizon
of a four-dimensional static black hole is a surface of constant curvature $E\in\{1,0,-1\}$ (for the metrics
$f_{\mathsf K}d\tau^2+f_{\mathsf K}^{-1}dr^2+r^2g_B$ of the proof of lemma \ref{lem:cross}, with $N=2$ and $\mathsf K=E$). A
sphere carries no TT tensor at all, so there is no perturbation of tensor type --- this is the classical fact that in
dimension four gravitational perturbations are of scalar and vector type only. On a flat or hyperbolic surface, every
TT tensor $h_B$ satisfies $\Delta_L^Bh_B=2Eh_B$: since $\Ric(g)=\frac12R(g)\,g$ for every metric on a surface, the
linearisation at $g_B$ in a TT direction gives $\frac12\Delta_L^Bh_B=Eh_B$. Hence $\lambda=2E$, and the radial
operator, given by the same computation as \eqref{radial}, reduces to $A_{2E}\psi=-r^{-N}(r^Nf_{\mathsf K}\psi')'$, which
is non-negative and has no normalisable zero mode. Degenerate Poincar\'e--Einstein metrics do exist in dimension four:
for \eqref{bh} with $N=2$, the period $\tau_0(r_h)=4\pi r_h/(1+3r_h^2)$ of $\tau$ is maximal at $r_h=1/\sqrt3$, where $\partial g_{r_h}/\partial r_h$ is a
non-zero $L^2$ element of $\ker(\Delta_L+2n)$ whose boundary variation vanishes. But this kernel is invariant under the
isometries of the conformal infinity, so it breaks no symmetry.

Pending a proof, a weaker statement is available: the extension holds at non-degenerate metrics
\cite[Prop.~A.1]{ACD2}, hence along any branch obtained by the implicit function theorem. The degree theory of
\cite{An3} cannot be substituted for the extension statement, since the computation of the degree in
\cite[Th.~B]{An3} uses \cite[Th.~2.5]{An3} to identify the fillings of the round $S^3$, as do the examples of
\cite[\S7]{An3}; that particular step is however independently available on spin manifolds, by the scalar curvature
rigidity theorem \cite[Th.~5.2]{AD}.

The results of \cite{ACD} which rest on \cite[Prop.~1]{ACD} (\cite[Prop.~2--6]{ACD} and \cite[Th.~3,~4]{ACD}) hold
as they stand for metrics on which the $S^1$ action is assumed rather than deduced. In the toral case, part of their
conclusions can be approached without any symmetry assumption through the Horowitz--Myers conjecture and its rigidity
case \cite{BH1,BH2}, although the translation is not literal. In the round case, where $\gamma=d\tau^2\oplus g_{+1}$
and $X_0=\partial_\tau$, it would suffice to prove that $\mathcal L_{\partial_\tau}g_{(3)}=0$, by remark
\ref{rem:weaker}\,(b), the manifolds considered in \cite{ACD} satisfying $\pi_1(M,\partial M)=0$.
\end{remark}

\subsection*{Acknowledgements}
The author thanks Michael Anderson for his kind response to a preliminary version of this paper.
This article was partially funded by the French National Research Agency ANR-23-CE40-0010: \emph{Einstein
constraints: past, present, and future} (EINSTEIN-PPF) and ANR-24-CE40-0702: \emph{Orbifold, scalar curvature and
rigidity} (OrbiScaR).

\subsection*{Data availability}
Data sharing is not applicable to this article, as no datasets were generated or analysed during the current study.

\subsection*{Conflict of interest}
The author has no conflict of interest to declare.

\subsection*{Statement on the use of artificial intelligence}
This article is the result of many exchanges between the author and an AI assistant (Claude, Anthropic), which was
used for the exploration and the elaboration of the proofs, for computations and numerical tests (the latter served
only to explore and to check: every statement of the paper is proved exactly, and no proof relies on a numerical
computation), for the search of references and for the writing of the text. The author takes full
responsibility for the content of the paper.

\appendix
\section{Computation of the modified potential}\label{app:pot}

In this appendix we give the details of the computation behind lemma \ref{lem:G}. We use the notation of section
\ref{sec:scalar}: $N=7$, $u=1/r$, $\mathsf x=r_h^{6}(1+r_h^2)u^{6}$, $\mathcal H=m+28\mathsf x$,
$\mathsf s=\mathsf xu^2$, $\hat F=fu^2=1+u^2-\mathsf s$, and $\mathsf y=m+\varpi$ as in section \ref{sec:scalar}.

\subsection*{A.1. The ground state transform}
Let $\Psi>0$ and $\omega=d\log\Psi/dr_*$, so that $\Psi^{-1}d^2\Psi/dr_*^2=d\omega/dr_*+\omega^2$. For $\Phi=\Psi\chi$,
$$
-\frac{d^2\Phi}{dr_*^2}+V_S\Phi=-\frac1\Psi\frac d{dr_*}\Big(\Psi^2\frac{d\chi}{dr_*}\Big)+\tilde V\Phi,\qquad
\tilde V=V_S-\frac{d\omega}{dr_*}-\omega^2 .
$$
If $\Phi$ solves \eqref{master}, multiplying by $\Phi=\Psi\chi$ and integrating by parts gives \eqref{gst}.

\subsection*{A.2. Derivatives}
Since $dr=-du/u^2$, $d/dr_*=f\,d/dr=-\hat F\,d/du$. Moreover
$$
\begin{gathered}
\hat F_u=2u-\frac{8\mathsf s}u,\qquad \mathcal H_u=\frac{6(\mathcal H-m)}u,\qquad \mathcal H_{uu}=\frac{30(\mathcal H-m)}{u^2},\\
V_S=\frac{\hat F\,U}{16\mathcal H^2},
\end{gathered}
$$
the last formula because $f/r^2=\hat F$. With $w^2=9+\mathsf y u^2$ we have $w_u=\mathsf y u/w$, and
$$
\frac d{du}\Big(w+\frac32\log\frac{w-3}{w+3}\Big)=\frac{\mathsf y u}w+\frac32\cdot\frac6{w^2-9}\cdot\frac{\mathsf y u}w
=\frac{\mathsf y u^2+9}{uw}=\frac wu ,
$$
so that $\mathsf D:=\frac d{du}\log\Psi_\varpi=\mathsf A-\frac{\mathcal H_u}{\mathcal H}$ with $\mathsf A=\frac{w-\frac12}u$. From
$\mathsf A_u=\frac{\mathsf y}w-\frac{w-\frac12}{u^2}$ and $\mathsf A^2=\frac{9+\mathsf y u^2-w+\frac14}{u^2}$,
$$
\mathsf A_u+\mathsf A^2=\mathsf y\Big(1+\frac1w\Big)+\frac{\frac{39}4-2w}{u^2} .
$$

\subsection*{A.3. The modified potential}
Then $\omega=-\hat F\mathsf D$, $d\omega/dr_*=\hat F(\hat F_u\mathsf D+\hat F\mathsf D_u)$, hence
$\tilde V=V_S-\hat F\hat F_u\mathsf D-\hat F^2(\mathsf D_u+\mathsf D^2)$,
with
$$
\mathsf D_u+\mathsf D^2=\mathsf A_u+\mathsf A^2-2\mathsf A\frac{\mathcal H_u}{\mathcal H}-\frac{\mathcal H_{uu}}{\mathcal H}+\frac{2\mathcal H_u^2}{\mathcal H^2}.
$$
This gives
$$
\tilde V=V_S-\hat F^2\Big(\mathsf A_u+\mathsf A^2-\frac{\mathcal H_{uu}}{\mathcal H}+\frac{2\mathcal H_u^2}{\mathcal H^2}\Big)
+\mathsf A\hat F\Big(\frac{2\hat F\mathcal H_u}{\mathcal H}-\hat F_u\Big)
+\frac{\hat F\hat F_u\mathcal H_u}{\mathcal H} .
$$
Multiplying by $16\mathcal H^2u^2w/\hat F$ and using $u\mathsf A=w-\frac12$ and the formula for $\mathsf A_u+\mathsf A^2$, we get
\begin{align*}
\frac{16\mathcal H^2u^2w}{\hat F}\tilde V={}&w\,u^2U-16\hat F\,\mathcal H^2\big(\mathsf y u^2(w+1)+w(\tfrac{39}4-2w)\big)\\
&-16\hat F\,w\,u^2\big(2\mathcal H_u^2-\mathcal H\mathcal H_{uu}\big)\\
&+16\,u\,w\big(w-\tfrac12\big)\big(2\hat F\mathcal H\mathcal H_u-\hat F_u\mathcal H^2\big)+16\,u^2w\,\hat F_u\mathcal H\mathcal H_u ,
\end{align*}
a polynomial in $u$, $w$, $m$ and $\mathsf x$. Substituting $\mathsf p=u^2$, $w=\vartheta+3$, $\mathsf s=\mathsf xu^2$ and $m=\vartheta(\vartheta+6)/\mathsf p-\varpi$ with $\varpi=5$, one
finds $\frac{16\mathcal H^2u^2w}{\hat F}\tilde V=u^{-4}\hat{\mathcal P}(\mathsf p,\vartheta,\mathsf s)$ with $\hat{\mathcal P}=\hat{\mathcal P}_0+\hat{\mathcal P}_1\mathsf s+\hat{\mathcal P}_2\mathsf s^2+\hat{\mathcal P}_3\mathsf s^3$,
\begin{align*}
\hat{\mathcal P}_0={}&576\vartheta^4+192\vartheta^5+16\vartheta^6\\
&+\mathsf p\,(3456\vartheta^2-12480\vartheta^3-10048\vartheta^4-2848\vartheta^5-352\vartheta^6-16\vartheta^7)\\
&+\mathsf p^2(-5760\vartheta+20560\vartheta^2+13120\vartheta^3+2560\vartheta^4+160\vartheta^5)\\
&+\mathsf p^3(2400-8800\vartheta-4000\vartheta^2-400\vartheta^3),\\
\hat{\mathcal P}_1={}&13824\vartheta^2+22656\vartheta^3+12032\vartheta^4+3552\vartheta^5+448\vartheta^6+16\vartheta^7\\
&+\mathsf p\,(-6912\vartheta-67072\vartheta^2-44864\vartheta^3-12480\vartheta^4-1056\vartheta^5)\\
&+\mathsf p^2(-3840+42400\vartheta+24320\vartheta^2+4880\vartheta^3),\\
\hat{\mathcal P}_2={}&64512\vartheta+98560\vartheta^2+53760\vartheta^3+14336\vartheta^4+896\vartheta^5\\
&+\mathsf p\,(-53760+25088\vartheta-19712\vartheta^2-17024\vartheta^3),\\
\hat{\mathcal P}_3={}&75264\vartheta+50176\vartheta^2+12544\vartheta^3 .
\end{align*}

\subsection*{A.4. The variable $\mathsf y$ and the factorisation}
By \eqref{tyu}, $\mathsf p=\vartheta(\vartheta+6)/\mathsf y$. Substituting and multiplying $\hat{\mathcal P}_i$ by $\mathsf y^{3-i}$, every coefficient of a power of $\vartheta$
becomes a polynomial in $\mathsf y$. For instance, the coefficients of $\vartheta^3$ and $\vartheta^4$ in $\mathsf y^3\hat{\mathcal P}_0$ are
$$
\begin{gathered}
20736\mathsf y^2-207360\mathsf y+518400=20736(\mathsf y-5)^2,\\
576\mathsf y^3-71424\mathsf y^2+671040\mathsf y-1641600=576(\mathsf y-5)^2(\mathsf y-114).
\end{gathered}
$$
In the same way, all the coefficients of $\mathsf y^3\hat{\mathcal P}_0$ vanish to order two at $\mathsf y=5$ and those of $\mathsf y^2\hat{\mathcal P}_1$ to order one,
and one obtains
\begin{align*}
\mathsf y^3\hat{\mathcal P}_0&=(\mathsf y-5)^2\vartheta^3\big[20736+(576\mathsf y-65664)\vartheta+(192\mathsf y-70848)\vartheta^2\\
&\qquad\qquad+(16\mathsf y-26976)\vartheta^3-4960\vartheta^4-448\vartheta^5-16\vartheta^6\big],\\
\mathsf y^2\hat{\mathcal P}_1&=(\mathsf y-5)\,\vartheta^2\big[(13824\mathsf y+27648)+(22656\mathsf y-296064)\vartheta\\
&\qquad\qquad+(12032\mathsf y-276096)\vartheta^2+(3552\mathsf y-101984)\vartheta^3\\
&\qquad\qquad+(448\mathsf y-16576)\vartheta^4+(16\mathsf y-976)\vartheta^5\big],\\
\mathsf y\,\hat{\mathcal P}_2&=\vartheta\big[(64512\mathsf y-322560)+(98560\mathsf y+96768)\vartheta+(53760\mathsf y-93184)\vartheta^2\\
&\qquad\qquad+(14336\mathsf y-121856)\vartheta^3+(896\mathsf y-17024)\vartheta^4\big],
\end{align*}
and $\hat{\mathcal P}_3=Q_3$. With $\mathsf y-5=m$ and $576\mathsf y+192\mathsf y\vartheta+16\mathsf y\vartheta^2=16\mathsf y(\vartheta+6)^2$, these are $m^2\vartheta^3Q_0$, $m\,\vartheta^2Q_1$ and $\vartheta Q_2$, and
$\mathsf y^3\hat{\mathcal P}=\mathcal P$, which proves lemma \ref{lem:G}.

\subsection*{A.5. Origin of the supersolution and the planar brane}
The supersolution is suggested by the asymptotic region: there, in the variable $u$ (asymptotic to the tortoise
coordinate $r_{*,\infty}-r_*$), $V_S=\frac{15}{4u^2}+m+O(1)$ (for the planar brane below, $V_S=\frac{15}{4u^2}+m+O(u^4)$).
With $q^2:=m+\varpi$, the model equation $-\Phi_{uu}+(\frac{15}{4u^2}+q^2)\Phi=0$ has the normalisable solution
$u^{1/2}I_2(qu)$, where $I_2$ and $I_3$ denote the modified Bessel functions of the first kind of orders $2$ and $3$,
and the classical bound
$I_3(z)/I_2(z)\ge z/(3+\sqrt{9+z^2})$ \cite{Amos} gives $\frac d{du}\log\big(u^{1/2}I_2(qu)\big)\ge\frac{w-1/2}u$; the factor
$1/\mathcal H$ accounts for the near-horizon region, and the shift $\varpi$, taken equal to $5$ in lemma \ref{lem:G}, for
the curvature of $B$.

For the planar black brane $f=r^2-r^{1-N}$, $N=7$, with momentum $q>0$ (flat horizon, $m=q^2$), the same supersolution with
$\varpi=0$ gives
$$
\begin{gathered}
\frac{16\mathcal H^2u^2}{\hat F}\tilde V=u^{4}\mathcal P_{\rm pl}+\frac{16\mathcal H^2u^2\hat F\,m(w-3)}{w(w+3)},\\
\mathcal P_{\rm pl}=\vartheta(\vartheta+6)\vartheta(768+576\vartheta+320\vartheta^2+16\vartheta^3)\\+\mathsf s\vartheta(21504+21504\vartheta+12544\vartheta^2+896\vartheta^3)+12544\,\mathsf s^2\vartheta(\vartheta+2),
\end{gathered}
$$
with $\hat F=1-u^8$, $\mathsf s=u^8$, $\vartheta=w-3$, so that $\vartheta(\vartheta+6)=mu^2$. All coefficients are non-negative, so the static scalar
master equation of the planar brane has no non-zero solution bounded at the horizon and normalisable at
infinity.

\section{Harmonics on an Einstein base}\label{app:harm}

The master equations of \cite{KI} hold for an Einstein base: this is observed in \cite{KI}, and \cite{KI2} treats
horizons which are Einstein manifolds from the outset. The Riemannian formulation of \cite{CDK}, the regularity and the
asymptotics of the master functions, and the reconstruction of a kernel element with vanishing master functions are
written there for a base of constant curvature; we justify here that they carry over to the Einstein base
$(B^N,g_B)$, $\Ric_B=Eg_B$, $E=N-1$. Let $D$ be the Levi-Civita connection of $g_B$; all
the operators below are those of $g_B$.

\begin{lemma}\label{lem:harm}
Let $(B^N,g_B)$ be a compact Einstein manifold, $\Ric_B=Eg_B$.

\textup{(i)} Let $S$ be a function with $\Delta S=\kappa S$, $\kappa>0$, and $(\nabla^2S)^\circ=\nabla^2S+\frac\kappa NSg_B$ the
trace-free part of its Hessian. Then
$$
\begin{gathered}
\delta(dS)=\kappa S,\qquad \Delta_HdS=\kappa\,dS,\qquad \nabla^*\nabla dS=(\kappa-E)dS,\\
\delta(Sg_B)=-dS,\qquad \delta(\nabla^2S)^\circ=\frac{N-1}N\Big(\kappa-\frac{NE}{N-1}\Big)dS,\\
\Delta_L(Sg_B)=\kappa\,Sg_B,\qquad \Delta_L(\nabla^2S)^\circ=\kappa\,(\nabla^2S)^\circ .
\end{gathered}
$$
Moreover $\|(\nabla^2S)^\circ\|^2=\frac{N-1}N\big(\kappa-\frac{NE}{N-1}\big)\kappa\|S\|^2$ in $L^2$.

\textup{(ii)} Let $V$ be a co-closed $1$-form with $\Delta_HV=\mu V$. Then $\delta^*V$ is trace-free and
$$
\nabla^*\nabla V=(\mu-E)V,\qquad \delta\delta^*V=\tfrac12(\mu-2E)V,\qquad \Delta_L\delta^*V=\mu\,\delta^*V .
$$
\end{lemma}
\begin{proof}
On an Einstein manifold $\Delta_Hd=d\Delta$, $\Delta_H=\nabla^*\nabla+E$ on $1$-forms, $\Delta_L(\varphi g)=(\Delta\varphi)g$,
$\Delta_L\delta^*=\delta^*\Delta_H$ \cite{Li61}, and $2\delta\delta^*\omega=\nabla^*\nabla\omega+d\delta\omega-E\omega$ (the Weitzenb\"ock
formula already used in the proof of lemma \ref{lem:simple}). Since $\nabla^2S=\delta^*dS$, this gives
$\Delta_L\nabla^2S=\kappa\nabla^2S$ and $\delta\nabla^2S=\frac12\big((\kappa-E)+\kappa-E\big)dS=(\kappa-E)dS$, hence
$\delta(\nabla^2S)^\circ=(\kappa-E-\frac\kappa N)dS$, which is the stated formula; the other identities of (i) are immediate,
and the norm is $\langle(\nabla^2S)^\circ,\nabla^2S\rangle=\langle\delta(\nabla^2S)^\circ,dS\rangle$ with $\|dS\|^2=\kappa\|S\|^2$. For (ii),
$\mathrm{tr}\,\delta^*V=-\delta V=0$ and the same formulas apply with $\delta V=0$.
\end{proof}

With $E=N-1$, the constant in (i) is $\frac{N-1}N(\kappa-N)=\frac{N-1}Nm$: the norm identity is the Lichnerowicz
inequality $\kappa\ge N$, with equality only if $(\nabla^2S)^\circ=0$, which characterises the round sphere \cite{Obata};
this is the exceptional case excluded in section \ref{sec:scalar}. In the notation of \cite{KI}, $k^2=\kappa$,
$\mathbb S_i=-k^{-1}D_iS$, $\mathbb S_{ij}=k^{-2}(\nabla^2S)^\circ$, and for vectors $k_V^2=\mu-E$,
$\mathbb V_{ij}=-k_V^{-1}\delta^*V$; lemma \ref{lem:harm} then gives the relations between harmonics used in
\cite{KI}, with the curvature $K=1$, for instance $D^j\mathbb S_{ij}=\frac{N-1}N\frac{k^2-N}k\mathbb S_i$ and
$D^j\mathbb V_{ij}=\frac{k_V^2-(N-1)}{2k_V}\mathbb V_i$, and $V$ is Killing if and only if $\mu=2E$.

\emph{How the base enters.} Let $g=g_2+r^2g_B$ be the warped product \eqref{bh}. By the formulas for the connection and
the curvature of a warped product, the curvature of $g$ is the sum of $r^2$ times the curvature of $g_B$, lifted to
$M$, and of terms which only involve $g_2$, $r$, $dr$, $\nabla^2r$ and $g_B$; in particular the components of the
curvature of $g$ with exactly one or exactly three indices on $B$ vanish, and those with two indices on $B$ and two on
the orbit space are $-r\,\nabla_a\nabla_br\,(g_B)_{ij}$. Consequently, for a symmetric $2$-tensor $h$ on $M$, the
curvature of $B$ enters $(\Delta_L+2n)h$ and $\beta(h)$ only through $\Ric_B=Eg_B$ and through the Lichnerowicz
Laplacian $\Delta_L^B=\nabla^*\nabla+2E-2\mathring R^B$ acting on the block of $h$ tangent to $B$; the other operators
of $B$ which appear are $\nabla^*\nabla$ on functions and on $1$-forms, $d$, $\delta$, $\delta^*$, the trace and the
multiplication by $g_B$. By lemma \ref{lem:harm}, all of them map the building blocks $S$, $dS$, $Sg_B$,
$(\nabla^2S)^\circ$ of the scalar type (respectively $V$, $\delta^*V$ of the vector type) to multiples of building blocks,
with coefficients depending only on $N$, $E$ and $\kappa$ (respectively $\mu$), which are those of a base of constant
curvature $1$. On a base of constant curvature, the rough Laplacian of $\mathbb S_{ij}$ or $\mathbb V_{ij}$ is used
instead of $\Delta_L^B$; on an Einstein base it is the latter which appears, and this is the only change
\cite[\S2]{KI2}. For the same reason, in the tensor sector, the eigenvalue of $\Delta_L^B$ replaces that of the rough
Laplacian, which gives \eqref{radial}. For a warped product with Einstein fibre over a \emph{one-dimensional} base, the
Lichnerowicz Laplacian of a TT tensor is computed in \cite[\S4]{DelayTT}; the base is here the two-dimensional orbit
space, so we argue directly.

Hence, for a perturbation of scalar or vector type, the linearised Einstein equations and the gauge condition reduce to
the same system of equations on the orbit space as for a base of constant curvature $1$, with the same
eigenvalue parameters. The gauge-invariant variables, the master functions and the master equations of
\cite[\S3, \S5]{KI}, the Riemannian formulation, the regularity and the asymptotics of the master function of an $L^2$
kernel element \cite[\S3]{CDK}, and the reconstruction of a kernel element with vanishing master function as
$\mathcal L_Wg$ with $W$ decaying \cite[\S4]{CDK}, are all obtained from this system by computations on the orbit
space; they therefore hold without change for the base $(B,g_B)$.

\end{document}